\documentclass[12pt]{article}
\usepackage{amssymb}
\usepackage{amsmath}
\usepackage{graphicx}
\usepackage{longtable}
\usepackage{psfrag}
\usepackage{mathrsfs}
\parskip \baselineskip
\newcommand{\ds}{\displaystyle}

\newcommand{\ol}{\overline}
\newcommand{\ul}{\underline}
\newcommand{\ra}{\rightarrow}

\newtheorem{thm}{Theorem}

\newtheorem{lm}{Lemma}

\newtheorem{prop}{Proposition}

\newtheorem{corol}{Corollary}

\newtheorem{ex}{Example}

\newtheorem{re}{Remark}

\begin{document}
\begin{center}
{\large \bf Lattices freely generated by posets within a variety.

Part II: Finitely generated varieties}

Jean Yves Semegni and Marcel Wild
\end{center}

\section{Introduction} \label{section1}
This article constitutes the second part of an essay dedicated to lattices freely generated by finite posets within a variety. The first part dealt with four "easy" cases, namely the variety of all semilattices, (general) lattices, distributive lattices, and Boolean lattices respectively. Special attention was paid to semilattices with a view to applications in Part II.

In the present Part II we are officially concerned  with finitely generated (f.g.) varieties $\cal{V}$ of lattices, in the usual sense that all subdirectly irreducible members are finite and, up to isomorphism, there are only finitely many of them. We wrote "officially" because quite a few preliminaries, hopefully interesting in their own right, will have to be digested before we come to f.g. varieties in section \ref{section6} and \ref{section7}.

The following problem is posed in section \ref{section2}: Given finite lattices $L_1, \cdots, L_t$, how much additional information about a subdirect product $L$ thereof is needed in order to compute $L$? Specifically, let $\pi_i : L\ra L_i$ be the $i$-th projection map  and $\sigma _i : L_i\ra L$ the corresponding smallest pre-image map. It turns out that the knowledge of the {\it connection maps} $\alpha_{i,j} : = \pi_i \circ \sigma _j$ from $L_j$ to $L_i\; (1\le i,j\le t)$ is sufficient, even if the $\pi_i$'s and $\sigma _i$'s themselves are unknown. Where one would get the connection maps from, will be seen in section \ref{section6}.

When computing $L$ from the connection maps it pays to replace the seemingly natural set $J(L)$ of join irreducibles by the larger {\it scaffolding} $G(L)$, which is defined as the union of the sets $\sigma_i(L_i\setminus\{0\})$. Following \cite{Wi2} we show in section \ref{section3}   that the $\vee$-semilattice freely generated by the partial $\vee$-semilattice $(G(L), \bigvee)$  is isomorphic to $L\setminus \{0\}$. The benefit is that free $\vee$-semilattices can be viewed as certain closure systems $\mathscr{C}$ which are amenable to the  implication $n$-algorithm introduced in Part I (it is fully discussed in \cite{mw1}). Namely, this algorithm is applicable whenever $\mathscr{C}$ is given by an implicational base $\Sigma$.

The scaffolding $G(L)$  contains the join core $K_{\vee}(L)$ which in turn contains $J(L)$. Section \ref{section4} investigates $K_{\vee}(L)$ when $L$ is a {\it modular} lattice. The view of $K_{\vee}(L)$ as linear hypergraph (pioneered in \cite{hw}) generalizes the projective geometry view of {\it complemented} modular lattices.
\vskip 0.25cm  
Section \ref{section5} fine-tunes the implication $n$-algorithm to the situation where  $\mathscr{C}$ is  isomorphic to a modular lattice, and thus consists of all order ideals of a poset which simultaneously are closed  with respect to some linear hypergraph. 
\vskip 0.25cm 
In section \ref{section6} finitely generated varieties $\mathcal{V}$  enter the stage. Let $L$ be the lattice freely generated within $\mathcal{V}$ by some finite poset $(P,\le)$. The calculation of $L$ is based on two essential ideas. 
\vskip 0.25cm 
First, $L$ is a subdirect product with factors $L_i$ from among the finitely many subdirect irreducibles of  $\mathcal{V}$. Crucially, since $L$ is free, the connection maps $\alpha_{i,j}$ between the $L_i$'s (section \ref{section2}) can be calculated in miraculous ways and they yield $G(L)$. 
\vskip 0.25cm 
Second, the fact that the partial semilattice $(G(L),\bigvee)$ freely generates $L$ as a semilattice, makes $L$ a closure system to which the  $(A,B)$-algorithm applies.

Section \ref{section7} focuses on the variety $\cal{V}$ which is generated by the smallest modular nondistributive lattice $M_3$.  For most of the $318$ posets $P$ with $\left|P\right|= 6 $ the cardinality of the $\cal{V}$-free lattice generated by $P$ was computed in \cite{yves}.

\section{Retrieving  a subdirect product from its connection maps} \label{section2}

Let $\phi : L \ra L_0$ be an epimorphism of lattices such that each $y \in L_0$ has a smallest pre-image $\sigma (y) = \bigwedge \{x \in L: \phi (x) = y\}$.  In particular that  is the case in our situation where all lattices are finite. It is shown in \cite{Wi2} that $\sigma : L_0 \ra L$ is an injective $\vee$-homomorphism.

Let $L\subseteq \ds \prod_{1\le i\le t}L_i$ be a subdirect product of lattices $L_1, \ldots, L_t$. Then all restricted projections $\pi _i : L \ra L_i$ have smallest pre-images $\sigma_i: L_i \ra L$ and 

\begin{eqnarray}\label{eq1}
 J(L) &=& \ds\bigcup_{1 \leq i \leq t} \sigma_i (J(L_i)).
\end{eqnarray}
This is implicit in \cite{Wi2} and explicitly in [1, Thm \;3.4].
Mutatis mutandis the same holds for meet irreducibles and biggest pre-images, but these will not concern us here. Observe that with $\sigma_j$ also $\alpha_{i,j}:=\pi_i\circ\sigma_j: L_j \ra L_i$ is a $\vee$-homomorphism. One readily checks that
\begin{eqnarray}\label{eq2}
 \alpha_{i,i}=id, \;\; \alpha_{i,j}\circ\alpha_{j,k}\le \alpha_{i,k}
\end{eqnarray}
for all $1\le i,j,k\le t$. Conversely, suppose one is given  lattices $L_1, \cdots, L_t$ and {\it any} family $\alpha_{i,j} : L_j \ra L_i$ of $\vee$-homomorphisms that satisfy (\ref{eq2}). Then these \textit{connection maps} $\alpha_{i,j}$ are induced by a suitable subdirect product as above. Namely, defining $\sigma_i': L_i \ra \ds \prod_{1\le i\le t}L_i $ by
\begin{eqnarray}\label{eq3}
\sigma'_i (y) &: = &(\alpha_{j,i} (y) | \ 1 \leq j \leq t)
\end{eqnarray}
the following takes place.
\begin{thm} \cite{Wi2} \label{thm1}
Given a set of connection maps satisfying (\ref{eq2}), there is a unique subdirect product $L\subseteq \ds \prod_{1\le i\le t}L_i$ such that the maps $\sigma'_i$ in (\ref{eq3}) are the smallest pre-image maps of the projections $\pi_i: L\ra L_i$ $(1\le i\le t)$.
\end{thm}
According to Theorem \ref{thm1} and (\ref{eq1}), the subdirect product $L$ can be calculated as the $\vee$-subsemilattice of $\ds \prod_{1\le i\le t}L_i $ generated by $\ds \bigcup_{1\le i\le t} \sigma'_i\left(J\left(L_i\right)\right)$. 

\begin{ex} \label{myexa}
{\rm Consider the two lattices} $L_1=\{a,b,\cdots, n\}$ {\rm and} $L_2=\{0,2,3,\cdots, 12,1\}$ {\rm (so $1$ is top) which are coupled by the} $\vee${\rm -homomorphisms} $\alpha_{1,2}$ {\rm and} $\alpha_{2,1}$ {\rm as indicated. For instance, as required in} (\ref{eq2}), {\rm we have}
$$(\alpha_{1,2}\circ \alpha_{2,1})(k) \;\;=\;\; \alpha_{1,2}(11)\;\;=\;\;c\;\;\le\;\; k. $$
\end{ex}

\begin{figure}[!h]
\begin{center}
\psfrag{b}{$b$}
\psfrag{c}{$c$}
\psfrag{d}{$d$}
\psfrag{f}{$f$}
\psfrag{h}{$h$}
\psfrag{i}{$i$}
\psfrag{e}{$e$}
\psfrag{l1}{$L_1\; =$}
\psfrag{a}{$a$}
\psfrag{g}{$g$}
\psfrag{j}{$j$}
\psfrag{k}{$k$}
\psfrag{m}{$m$}
\psfrag{alpha2,1}{$\alpha_{2,1}$}
\psfrag{alpha1,2}{$\alpha_{1,2}$}
\psfrag{lambda1}{$\Lambda_1:$}
\psfrag{lambda2}{$\Lambda_2:$}
\psfrag{l2}{$L_2$} \psfrag{l1}{$L_1$}
\includegraphics[scale=.45]{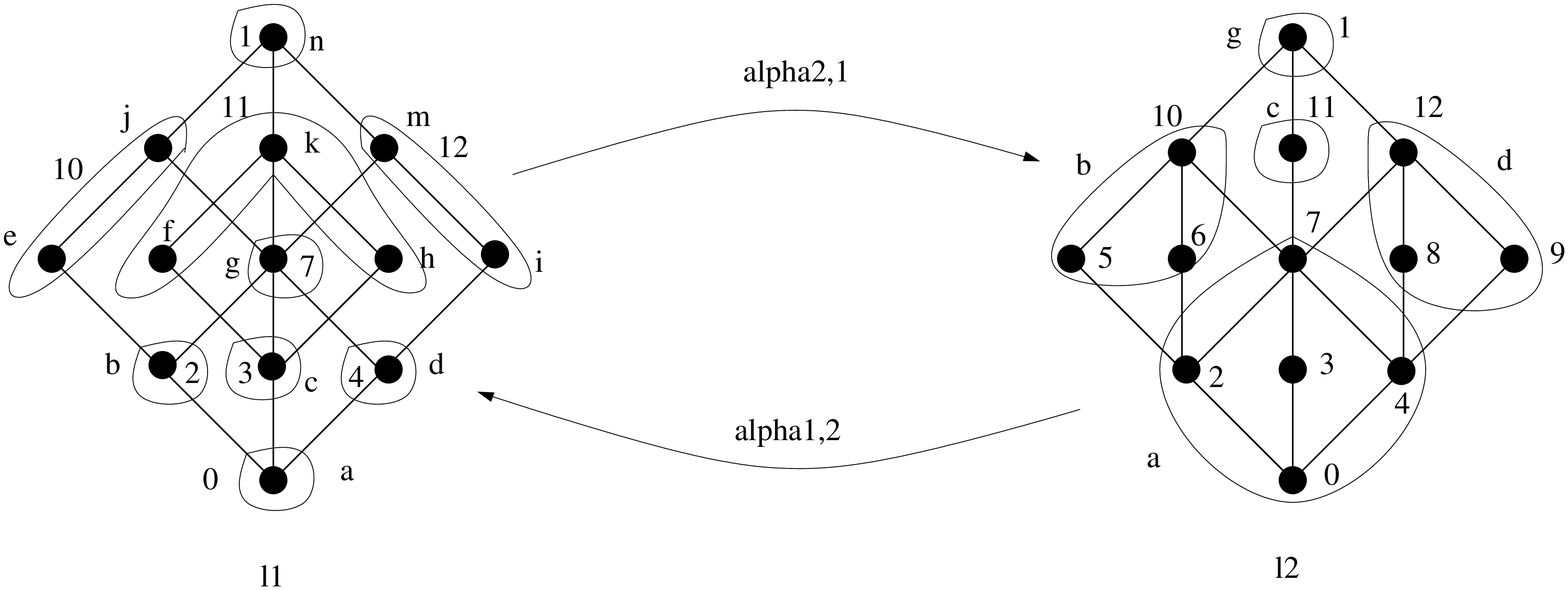}
\caption{}
\label{fig1}
\end{center}
\end{figure}

One verifies that:
\begin{eqnarray} \label{eq4}
J\left(L_1\right)&=&\{b, c, d, e, f, h, i\} \nonumber \\
J\left(L_2\right)&=&\{2, 3, 4, 5, 6, 8, 9, 11\}\nonumber \\
						\nonumber\\
\sigma'_1(b) \;\;=\;\;(b,2)  && \sigma'_2(2)\;\;=\;\;(a,2)\nonumber \\
\sigma'_1(c) \;\;=\;\;(c,3)  && \sigma'_2(3)\;\;=\;\;(a,3)\nonumber \\
\sigma'_1(d) \;\;=\;\;(d,4)  &&\sigma'_2(4)\;\;=\;\;(a,4) \nonumber \\
\sigma'_1(e) \;\;=\;\;(e,10)  &&\sigma'_2(5)\;\;=\;\;(b,5) \\
\sigma'_1(f) \;\,=\;\;(f,11)  &&\sigma'_2(6)\;\;=\;\;(b,6)\nonumber \\
\sigma'_1(h) \;\,=\;\;(h,11)  &&\sigma'_2(8)\;\;=\;\;(d,8)\nonumber \\
\sigma'_1(i) \;\,\;=\;\;(i,12)  &&\sigma'_2(9)\;\;=\;\;(d,9)\nonumber \\
                             &&\sigma'_2(11)\,=\;\,(c,11)\nonumber
\end{eqnarray}

It turns out that here taking suprema of {\it pairs} (as opposed to triplets, quadruplets, $\cdots$) of elements of (\ref{eq4}) suffices to generate the subdirect product $L\subseteq L_1\times L_2$. 
\newpage
\begin{figure}[!h]
\begin{center}
\includegraphics[scale=.375]{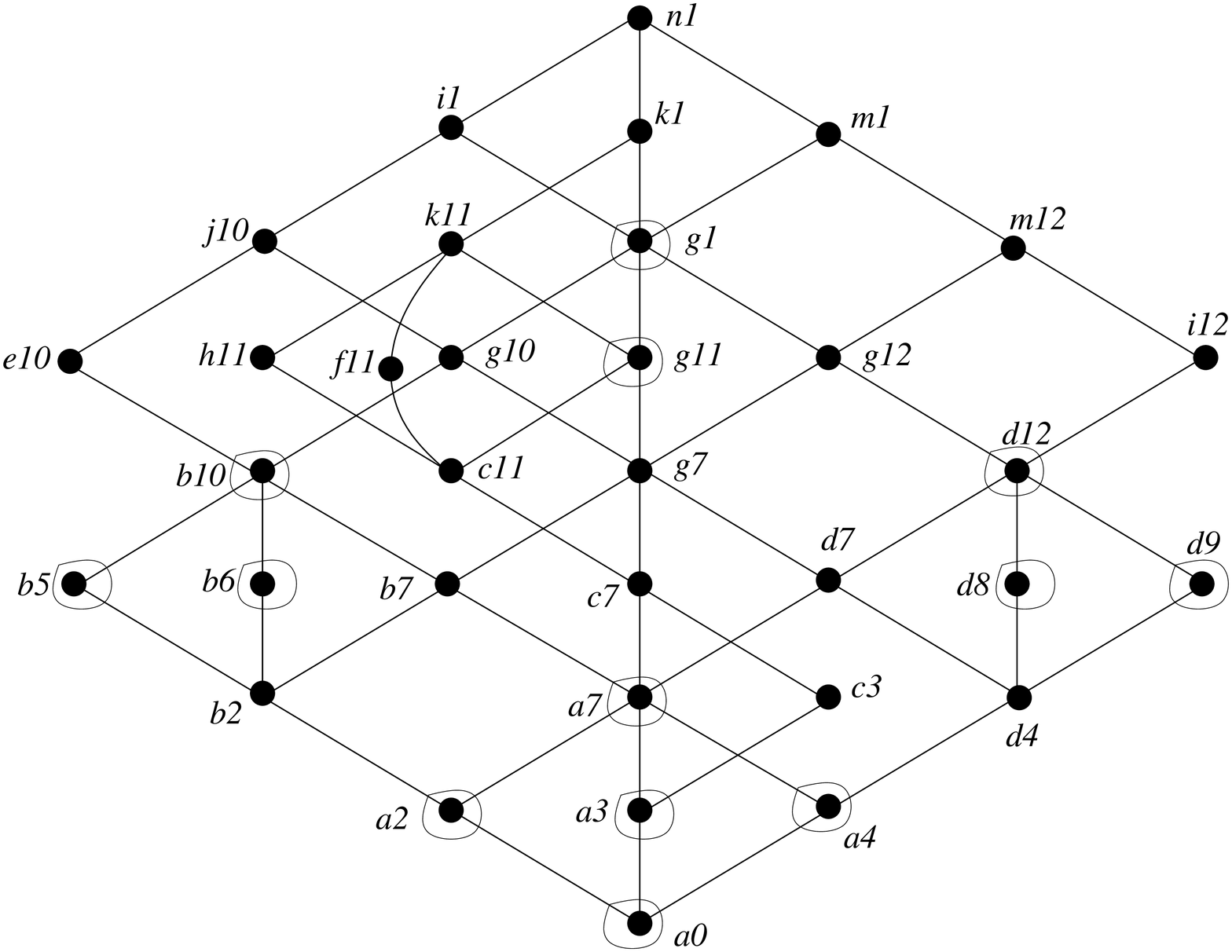}
\caption{}
\label{fig2}
\end{center}
\end{figure}

When generating $L$ from $J(L)$, some elements of $L$, say $(g,1)$, may be duplicated many times: 

\hspace*{2cm}$(b,5)\vee (d,8)\;\;\;\,=\;\;(b,5)\vee (d,9)\;=\;(b,6)\vee (d,8)\,\;\;=\;(b,6)\vee (d,9)$

\hspace*{1.75cm}=\;$(c,11)\vee(b,5)\;=\;(c,11)\vee (b,6)\;=\; (c,11)\vee (d,8)\;=\; (c,11)\vee (d,9)$
=\; $(g,1).$

This  is inefficient and a better way will be approached in section \ref{section3}. 

\section{The algorithmic advantage of the scaffolding $G(L)$ over $J(L)$} \label{section3}

It is straightforward to design a general purpose algorithm for calculating the subalgebra generated by a given subset  of a universal algebra. Using hashing techniques the inefficient regeneration of elements, as in Example \ref{myexa}, can be partly cured. For the case where the universal algebra is a not too large lattice, this approach has been taken (among other methods) in \cite{berman}.  Similarly the authors of \cite{lux} proceed to compute the lattice $L$ of all submodules of a module. We note that some of the theory developed in \cite{hw} is rediscovered. 

In contrast, our philosophy is the following. For $a\in L$ put $J(a) := \{p\in J(L)\;|\; p\le a\}$. We shall identify $L$ with the isomorphic closure system $\mathscr{C} =\{J(x)\;|\;x\in L\}$ and seek some suitable implicational base $\Sigma$ for $\mathscr{C}$. The point is that  with the implication $n$-algorithm of Part I, $\mathscr{C}$ can be computed faster  as the set $\mathscr{C}\left(\Sigma\right)$ of all $\Sigma$-closed subsets of $J(L)$. 

One way to come up with such a $\Sigma$ is as follows.   For any lattice $L$, let $R\subseteq L$
be such that it contains the join core $K_{\vee}(L)$. Recall from Part I, section $5$, that an implicational base $\Sigma$ of $\mathscr{C}\cong L$ is then obtained by collecting all the implications $A\ra J\left(\bigvee A\right)$ where $A$ ranges over those subsets of $J(L)$ for which $\bigvee A \in R$. In particular, $\Sigma$ contains all the implications $\{p\}\ra J(p)$ with $p$ ranging over $J(L)$. 

In this section we shall exhibit a convenient set $R=G(L)$ by merely exploiting that $L$ is subdirectly {\it reducible}. Thus let $L$ be a subdirect product of lattices $L_1, \cdots, L_t$ where, additionally to section \ref{section2}, the $L_i$'s must be subdirectly irreducible.

Akin to (\ref{eq1}) we define the {\it scaffolding} ("Ger\"{u}st" in \cite{Wi2})  of $L$ as
\begin{eqnarray}\label{eq5}
 G(L) &: = & \ds\bigcup_{1 \leq i \leq t} \sigma_i \left(L_i \setminus \{\emptyset\}\right).
\end{eqnarray}
Despite appearances, $G(L)$ is not dependent on the particular subdirect decomposition of $L$. For modular $L$ this will be shown  in section \ref{section4}.  As  seen in Part I,  as a subset  of $L$  it automatically becomes a partial semilattice $(G(L), \bigvee)$. It turns out that  $(G(L), \bigvee)$ {\it freely} generates $L$, that is,
\begin{eqnarray}\label{eq6}
F_{\vee}\left(G(L), \bigvee \right) &\simeq & L\setminus \{\emptyset\}.
\end{eqnarray}
Here comes the proof of (\ref{eq6}), which essentially is a translation of [3, 1.3]: 

For all $x \in L$ the full ideal $\varepsilon (x): = \{a \in G(L) : a \leq x \}$ contains $J(x)$ because of (\ref{eq1}) and (\ref{eq5}). Therefore the map $\varepsilon: L \ra F_{\vee}\left(G(L),\bigvee)\right)$ satisfies
\begin{eqnarray*}
\left(\forall x,z \in L\right) \;\;\; \left(x\le z \iff \varepsilon(x)\subseteq \varepsilon (z)\right).
\end{eqnarray*}
 In order to show that $\varepsilon $ is {\it onto} and hence an isomorphism of semilattices, consider
 $\kappa_i(x) : = \sigma_i(\phi_i (x))$. Then $\kappa_i : L \ra L$ is a kernel operator, i.e. is anti-extensive, idempotent, and preserves suprema. Fix $A \in F_{\vee}\left(G(L), \bigvee \right)$ and any $a$ in $\varepsilon (\bigvee A) \supseteq A$. We need to show that $a\in A$. For some $i \in \{1, \cdots, t\}$, we have $a = \kappa_i (a) \leq \kappa_i (\bigvee A) = \bigvee \kappa_i (A)$. Thus $\bigvee \kappa_i(A)\in G(L)$, and $\kappa_i (A) \subseteq A$ since $A$ is hereditary. Therefore $\ds\bigvee \kappa_i (A) \in A$ by definition of  $\bigvee$-ideal. Because $A$ is  hereditary,  $a \in A$.

\section{The join core of  modular lattices} \label{section4}

 As in Part I the natural closure operator associated to a lattice $L$ maps $X\subseteq J(L)$ to $\ol{X}:=J\left(\bigvee X\right)$. Recall that the set $E(L)$ of {\it essential} elements consists of those $x\in L$ which (alias $J(x)$) contain a proper quasiclosed generating set. Further the {\it join core} $K_{\vee}(L)$ is $E(L)$\begin{tabular}{l}  \includegraphics[scale =.4]{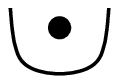} \end{tabular}$A(L)$ with $A(L)\subseteq J(L)$ the set of atoms of $L$. From (\ref{eq6}) and Theorem $2$ in Part I follows that the scaffolding comprises  $K_{\vee}(L)$. In fact,  according to [9, Thm. 7] one has 
\begin{eqnarray}\label{eq7}
 \quad K_\vee (L) &=& \ds\bigcup_{1 \leq i \leq t} \sigma_i \left(K_\vee (L_i)\right)
\end{eqnarray}
 where $L$ is any subdirect product of  subdirectly irreducible lattices $L_i \;(1\le i\le t)$. The proof of (\ref{eq7}) involves an application of Duquenne's multi-purpose $M_3$-$N_5$-lemma.

Here we show that for modular lattices  $K_{\vee}(L)$ is readily found. Along the way a proof of (\ref{eq7}) in the modular case unfolds. In fact the union in  (\ref{eq7}) turns out to be disjoint. Lemma \ref{lm1} and  Theorem \ref{thm4} below are based on [3, p.58];  we use the opportunity to mend some minor typos and expand some arguments.

\begin{lm} \label{lm1}  \cite{Wi2}:
Let $\phi : L \ra L_0$ be a lattice epimorphism with smallest pre-images $\sigma : L_0 \ra L$, and fix a nonzero element of type $v=\sigma\phi(v)$. Let $v/w$ be any prime quotient and let $r/s$ be a prime (i.e. covering) quotient in $L_0$ which is projective to the prime quotient $\phi v/\phi w$. Then $\ul{r}=\sigma r$ and $\ul{s}=\bigvee \{u\le \ul{r}\,|\,\phi u \le s\}$ yield a prime quotient $\ul{r}/\ul{s}$ which is projective to $v/w$.
\end{lm}
{\bf Proof:} Whereas $\ul{r}=\sigma r$, generally $\ul{s}\ne \sigma s$.  Rather $\ul{s}$ is the greatest element below $r$ that maps to $s$. In particular $\ul{r}/\ul{s}$ is a {\it prime} quotient (similarly for other quotients to come). One readily verifies that with $v/w$ also $\phi v/\phi w$ is a prime quotient. By assumption there are prime quotients $v_i/w_i \;(0\le i\le n)$ in $L_0$ such that $v_i/w_i$ is transposed to $v_{i-1}/w_{i-1} \;(1\le i\le n)$ and such that $v_0/w_0=r/s$ and $v_n/w_n = \phi v/\phi w$.  For all $0\le i\le n$ we put $\ul{v}_i=\sigma v_i$ and $\ul{w}_i:=\bigvee \{u\in L | u\le \ul{v}_i, \phi u\le w_i\}$. Notice that $\ul{v}_n/\ul{w}_n = v/w$ since $w$ is obviously the largest element below $v$ with $\phi$-image $\le w_n$. Ditto, $\ul{v}_0/\ul{w}_0=\ul{r}/\ul{s}.$ 

In order to see that $\ul{v}_i/\ul{w}_i$ is  transposed to $\ul{v}_{i-1}/\ul{w}_{i-1}\; (1\le i\le n)$, assume w.l.o.g. that $v_{i-1}/w_{i-1}$ transposes {\it up} to $v_i/w_i$, that is, $v_{i-1}\wedge w_{i}=w_{i-1}$ and $v_{i-1}\vee w_i=v_i$. We conclude that $\phi(\ul{v}_{i-1}\wedge \ul{w}_i) =\phi (\ul{v}_{i-1}) \wedge \phi (\ul{w}_i) = v_{i-1}\wedge w_i=w_{i-1}$.
Hence $\ul{v}_{i-1}\wedge \ul{w}_i\le \ul{w}_{i-1}$ by definition of $\ul{w}_{i-1}$. The  inequality $\ge$ follows from $\ul{w}_{i-1}\le \ul{v}_{i-1}$ and (clearly) $\ul{w}_{i-1}\le \ul{w}_i$. Recalling that $\sigma$ is a join homomorphism it follows from $v_{i-1}\vee w_i=v_i$ that $\ul{v}_{i-1}\vee \sigma w_i = \ul{v}_i$. Hence $\ul{v}_{i-1}\vee \ul{w}_i=\ul{v}_i$. \hfill $\blacksquare$

Call an element $v\ne 0$ of a lattice $L$ {\it sub-irreducible} if all prime quotients  $v/w$ are mutually projective. For instance each join irreducible is sub-irreducible.
\begin{thm}\label{thm4}\cite{Wi2}:
Let $L$ be a modular lattice. Then $G(L)$ consists of all sub-irreducible elements. The lattices $L_i \;(1\le i\le t)$ in (\ref{eq5}) are unique up to isomorphism and the union in (\ref{eq5}) is disjoint.
\end{thm}
{\bf Proof:} Let us fix any sub-irreducible $x\in L$ and argue why we must have $x\in \sigma_i(L_i\setminus\{0\})$ for some $i$. Assume to the contrary that $x\notin \sigma_i(L_i\setminus\{0\})$ for all $i$. Then $x$ is not the smallest element of any $\theta_i$-class (where $\theta_i$ is the kernel of $L \ra L_i$), and so for each $i$ there is some lower cover $y \prec x$  which is in the same $\theta_i$-class as $x$. It cannot be that for {\it all} $i$ {\it all} $y \prec x$ are in the same $\theta_i$-class as $x$  because then $\bigwedge \{\theta_i| \ 1 \leq i \leq t \} \neq 0$, contradicting the fact that $L$ is a subdirect product of the $L_i$'s. Hence there is an $i$ and lower covers $y_1, y_2$ of $x$ such that $(x, y_1) \in \theta_i$ but $(x,y_2) \not\in \theta_i$. Yet this cannot be since $x/y_1$ by assumption is projective to $x/y_2$. It follows that $\{x \in L| \ x$ is {\it sub-irreducible}$\}$ is a subset of $G(L)$.  So far, modularity was not used.

Conversely, pick $v\in G(L)$, say $v=\sigma_i \phi(v)$. Let $v/w$ be a prime quotient and let $\phi v/\phi w$ be its image in $L_i$. Since $L_i$ is modular, it is simple\footnote{In the nonmodular case $L_i$ is merely subdirectly irreducible. Finding $r/s$ becomes more subtle and also {\it weak} projectivities must be dealt with [3, p.58].}, and so {\it any} fixed prime quotient $r/s$ in $L_i$ will be projective to $\phi v/\phi w$. By Lemma \ref{lm1} there is a prime quotient $\ul{r}/\ul{s}$ that is projective to $v/w$. Crucially, since $\ul{r}/\ul{s}$ depends on $r/s$ (and not on $v/w$), $\ul{r}/\ul{s}$ is projective to any other prime quotient $v/w'$ as well. Hence $v$ is sub-irreducible. 

It is well known  that $L$ being modular the projectivity classes of prime quotients correspond bijectively to the simple factors $L_i \; (1\le i\le t)$. This establishes the disjointness of the union in (\ref{eq5}). \hfill $\blacksquare$

\begin{figure}[!h]
\begin{center}
\psfrag{b}{$b$} \psfrag{c}{$c$} \psfrag{d}{$d$} \psfrag{f}{$f$} \psfrag{e}{$e$} \psfrag{i}{$i$} \psfrag{h}{$h$} \psfrag{b2}{$(b,2)$} \psfrag{c3}{$(c,3)$} \psfrag{d4}{$(d,4)$} \psfrag{f11}{$(f,11)$} \psfrag{h11}{$(h,11)$} \psfrag{i12}{$(i,12)$} \psfrag{e10}{$(e,10)$} \psfrag{a4}{$(a,4)$}  \psfrag{b6}{$(b,6)$} \psfrag{b5}{$(b,5)$} \psfrag{c11}{$(c,11)$}  \psfrag{d9}{$(d,9)$} \psfrag{a2}{$(a,2)$} \psfrag{d8}{$(d,8)$} \psfrag{a3}{$(a,3)$} \psfrag{lambda}{$\Lambda \;:$} \psfrag{lambdaprime1}{$\Lambda'_1 \; :$} \psfrag{lambdaprime2}{$\Lambda'_2 \; :$}
\includegraphics[scale=.35]{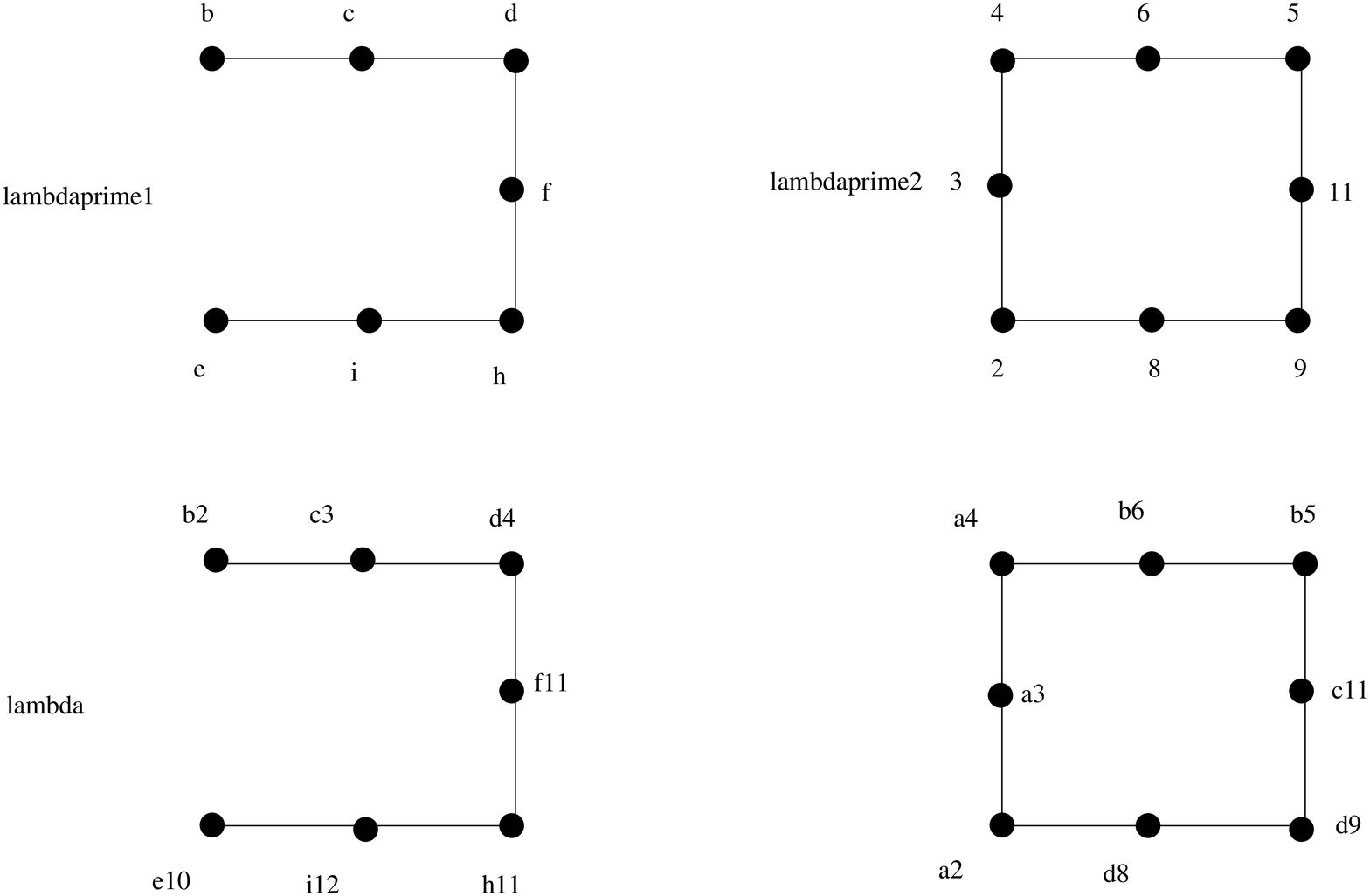}
\caption{}
\label{fig3}
\end{center}
\end{figure}

In order to succinctly describe the subset $ K_{\vee}(L)$ of $G(L)$ in the modular case, call  $x\in L$ a {\it line top}\footnote{This crisp name has been recently introduced by C. Herrmann; in \cite{hw} and elsewhere line tops were called $M_n$-elements.}  if $x$ has  $n \ge 3$ lower covers $x_i$ and their meet $\ul{x}$ is a lower cover of each $x_i$. So far $L$ could be any finite lattice. Since  all prime quotients of the interval sublattice $[\ul{x}, x]$ are clearly mutually projective,  each line top $x$ is in $G(L)$. Actually nonclosed quasiclosed generating sets of $J(x)$ are easy to find [13, p.156], and so even $x\in E(L)$. For modular lattices $L$ the line tops are the \textit{only} reducible essential elements. This was first shown in [9, Thm.9], other proofs are mentioned in [13, p.157]. Hence it follows from $K_{\vee}(L)=E(L)$\begin{tabular}{l}\includegraphics[scale =.4]{union.eps}\end{tabular}$A(L)$ that 
\begin{eqnarray}\label{eq8}
K_{\vee}(L) &=& J(L) \begin{tabular}{l} \includegraphics[scale =.4]{union.eps} \end{tabular}\{x\in L\;|\; x\; \,{\rm is\; \, a\; \, line \; \, top}\} 
\end{eqnarray}

for each finite modular lattice. As in (\ref{eq1}) the line tops occuring in any of the simple factors $L_i$ match the line tops of $L$,  and by Theorem \ref{thm4} the union in (\ref{eq7}) is disjoint. An at least $3$-element subset $l\subseteq J(L)$ maximal with the property that $p\vee q =\bigvee l$ for all distinct $p,q \in l$, is called a {\it line} of $L$. It is shown in \cite{hw} that  the line tops $x$ of $L$ are exactly the elements of type $x=\bigvee l$ with $l$ a line. One can have $x=\bigvee l_1 = \bigvee l_2$ for $l_1\ne l_2$. Furthermore $|l_1\cap l_2| \le 1$ for all lines $l_1\ne l_2$. A collection $\Lambda$ of lines for which  each line top $x$ contains exactly one $l\in \Lambda$ with $\bigvee l = x$, is called a {\it base of lines}. 
\begin{ex}\label{ex2}
{\rm Consider the lattice $L$ in Example \ref{myexa} which happens to be modular and which is a subdirect product of the simple lattices $L_1$, $L_2$ in Fig.\ref{fig1}.  The line tops of $L_1$ are $g, k, n$. The only line for $g$ is $l=l(g)=\{b,c,d\}$. The line top $k$ houses the two lines $\{f, b, h\}$ and $\{f,d,h\}$. Let us pick, say, $l(k)=\{f,d,h\}$. Similarly, say $l(n)=\{e, h, i\}$. The resulting base of lines is $\Lambda _1'=\{l(g), l(k), l(n)\}$. In the same way, one possible base of lines for $L_2$ is $\Lambda _2'=\{l(7), l(10), l(12), l(1)\}$ (e.g. $l(10) =\{4,5,6\}$). They are shown in Fig.\ref{fig3}. The corresponding base of lines for $L$ is $\Lambda: = \Lambda _1 \cup \Lambda _2$, where $\Lambda _i := \sigma _i(\Lambda' _i)$ is defined in the obvious way (see Fig.\ref{fig3}). Generally the number of connected components of any base of lines of a modular lattice equals the number of its simple factors. }
\end{ex} 
\newpage
For instance, the pre-image of the line top $12 \in L_2$, which can be evaluated as
\[
\sigma _2(12)\;=\;\sigma _2(2\vee 9) \;=\;\sigma _2(2)\vee \sigma _2(9)\;=\; (a,2)\vee (d,9)\;=\;(d,12),
\]
is a line top of $L$ and hence in $K_{\vee}(L)$. From Fig.\ref{fig2} one sees (after a while) that the element $(j,10)$ is sub-irreducible, i.e. the quotient $(j,10)/(e,10)$ is projective to $(j,10)/(g,10)$. According to Theorem \ref{thm4}  this forces $(j,10)\in G(L).$ Indeed, one checks that $(j,10)\;=\;\sigma _1(j)\in \sigma _1\left(L_1\setminus \{0\}\right)$. Clearly $(j,10)\notin K_{\vee}(L)$ since $(j,10)$ is neither join-irreducible nor a line top. Notice that e.g. $(b,7)\in L\setminus G(L)$.

\section{Algorithmic details in the modular case}\label{section5}
Let $L$ be any modular lattice with a base of lines $\Lambda$.
As seen in section \ref{section4}, the join core $R=K_{\vee}(L)$ consists of all the join irreducibles and the line tops $\overline{l}$ corresponding to the lines $l\in \Lambda$. Therefore (section \ref{section3}) an implicational base $\Sigma$ of the closure system $\mathscr{C} = \{J(x)\;|\; x\in L\}$  isomorphic to $L$ is obtained by taking all implications $A\ra J(\bigvee A)$ where $A$ is such that $\bigvee A \in R$. Here, besides $\{p\}\ra J(p) \; (p\in J(L))$, it suffices to take the implications $A\ra J(\bigvee A)$ of type $l\ra J(\overline{l})$ with $l\in \Lambda$. The purpose of  section \ref{section5} is to exploit this special type of $\Sigma$ in order to speed up the implication $n$-algorithm presented in Part I. The section is quite technical and may be skipped without loss of continuity.

To fix ideas, let us return to the lattice $L$ of Example \ref{ex2}. Put $J_i :=\sigma _i\left(J(L_i)\right)$ for $i=1,2$. The family $J[\Lambda]$ of all $\Lambda$-closed subsets $Z\subseteq J$ consists exactly of the sets $Z=X$\begin{tabular}{l} \includegraphics[scale =.4]{union.eps} \end{tabular}$Y$ where $X$ and $Y$ range over the accordingly defined families $J_1[\Lambda _1]$, respectively $J_2[\Lambda _2]$. Hence it makes sense to determine $J_1[\Lambda _1]$ and $J_2[\Lambda _2]$ apart, and {\it afterwards} worry to weave in all implications ${p}\ra J(p)$.

We identify subsets of $J_1 \, = \; \{(b,2),\cdots, (e,10)\}$ with their characteristic $0,1$-vectors but besides $0,1$ introduce other
symbols $2, l, \delta, \varepsilon$ in order to get {\it multi-valued rows} that compactly encode certain families of subsets of $J_1$. For simplicity we write $b2$  for $(b,2)$ and so forth. For starters, let $l_1=\{b2,c3,d4\}$ be our first line from  $\Lambda _1=\{l_1,l_2,l_3\}$. The family $J_1[l_1]\,:= \;\left\{X\subseteq J_1\;:\;\left|X\cap l_1\right|\in {0,1,3}\right\}$ of all $l_1$-closed subsets can be represented by the first (multivalued) row below in this (multivalued) context:

\begin{table}[!h]
\begin{center}
\begin{tabular}{|l|l|l|l|l|l|l|} 
$b2$ & $c3$ & $d4$ & $f11$ & $h11$ & $i12$ & $e10$ \\
\hline 
$ l$ & $l $  & $l $ & $2 $ & $2 $ & $ 2$ & $ 2$\\
\hline 
$ \varepsilon$ & $\varepsilon $  & $ {\bf 0}$ & $2 $ & $2 $ & $2 $ &  $2$\\
\hline  
$ \delta $ & $\delta $  & ${\bf 1} $ & $2 $ & $2 $ & $2 $ & $2$\\
\hline   
 &   &  & &  &  & \\
\hline   
$ \varepsilon _1$ & $\varepsilon _1 $  & $0 $ & $\varepsilon _2 $ & $\varepsilon _2 $ & $2 $ & $2$ \\
\hline   
$ \delta _1$ & $ \delta _1$  & $1 $ & $ \delta _2$ & $ \delta_2$ & $2 $ &  $2$\\
\hline   
&   &  &  &  &  & \\
\hline   
$ \varepsilon _1$ & $\varepsilon _1 $  & $ 0$ & $2 $ & $ {\bf 0}$ & $\varepsilon _2 $ & $\varepsilon _2 $\\
\hline   
$\varepsilon  $ & $ \varepsilon $  & $0 $ & $0 $ & $ {\bf 1}$ & $ \delta $ & $\delta$\\
\hline   
$ \delta $ & $\delta $  & $1 $ & $0 $ & ${\bf 0} $ & $\varepsilon$  & $\varepsilon $\\
\hline   
$ \delta _1$ & $\delta _1 $  & $1 $ & $1 $ & ${\bf 1} $ & $\delta _2$  & $\delta _2 $\\
\hline  
\end{tabular}
\caption{}
\label{tab1}
\end{center}
\end{table}

When we identify subsets with their characteristic vectors, then $lll$ in the first row is a shorthand for the family $\{000, 100, 010, 001, 111\}$ of all $l_1$-closed subsets of $\{b2,c3,d4\}$, and a symbol $2$ at any position means that the corresponding element is free to be present or not. We need two more symbols.\newline

Let $\varepsilon \varepsilon$ be a shorthand for $\{00, 01, 10\}$ ("at most one $1$") and $\delta \delta$ a shorthand for $\{00, 11\}$ ("dichotomy": all $1$'s or all $0$'s). The second and third row encode the fact that $J_1[l_1]=\mathscr{F}_1$\begin{tabular}{l} \includegraphics[scale =.4]{union.eps} \end{tabular}$\mathscr{F}_2 $  where 
\begin{eqnarray*}
\mathscr{F}_1  &:=&\left\{X\in J_1[l_1]:\, d4\notin X\right\} \;\;=\;\;\left\{X\subseteq J_1:\, \left|X\cap \{b2, c3\}\right|\in\{0,1\}\right\} \\
\mathscr{F}_2  &:=&\left\{X\in J_1[l_1]:\, d4\in X\right\} \;\;=\;\;\left\{X\subseteq J_1:\, \left|X\cap \{b2, c3\}\right|\in\{0,2\}\right\}
\end{eqnarray*}
 Consider the next line $l_2=\{d4, f11, h11\}$ of $\Lambda _1$. The fact that we did split $\mathscr{F}$ with respect to $d4$, which is the intersection of $l_1$ and $l_2$, benefits the imposition of $l_2$. Thus the fourth and fifth row encode all $X\subseteq J_1$ which are $\{l_1, l_2\}$-closed. Splitting each row with respect to $h11$ (the intersection of $l_2$ and $l_3 =\{h11, i12, e10\}$) yields the last four rows. They encode the family $J_1[\Lambda _1]$ of all $\Lambda _1$-closed sets $X\subseteq J_1$.

As seen, since $\sigma$ is injective and join-preserving, the poset $(J_1,\le)$ induced  by $(L,\le)$ is isomorphic to the poset $\left(J(L_1),\le \right)$ induced  by $(L_1,\le)$. Ditto $(J_2,\le)\simeq \left(J_2(L_2),\le\right)$. However, $\left(J(L),\le \right)$ features {\it more} comparabilities than the disjoint union of $(J_1,\le)$ and $(J_2,\le)$: 

\newpage
\begin{figure}[!h]
\begin{center}
\psfrag{e10}{$e10$} \psfrag{f11}{$f11$} \psfrag{h11}{$h11$} \psfrag{i12}{$i12$} \psfrag{b2}{$b2$} \psfrag{c3}{$c3$} \psfrag{d4}{$d4$} \psfrag{b5}{$b5$} \psfrag{b6}{$b6$} \psfrag{c11}{$c11$} \psfrag{d8}{$d8$} \psfrag{d9}{$d9$} \psfrag{a2}{$a2$} \psfrag{a3}{$a3$} \psfrag{a4}{$a4$} \psfrag{J}{$\left(J(L),\le\right)\;\;= \;\;$} \psfrag{J1}{$\left(J_1,\le\right) \;\;=$} \psfrag{J2}{$\left(J_2,\le\right) \;\;=$} 
\includegraphics[scale=.4]{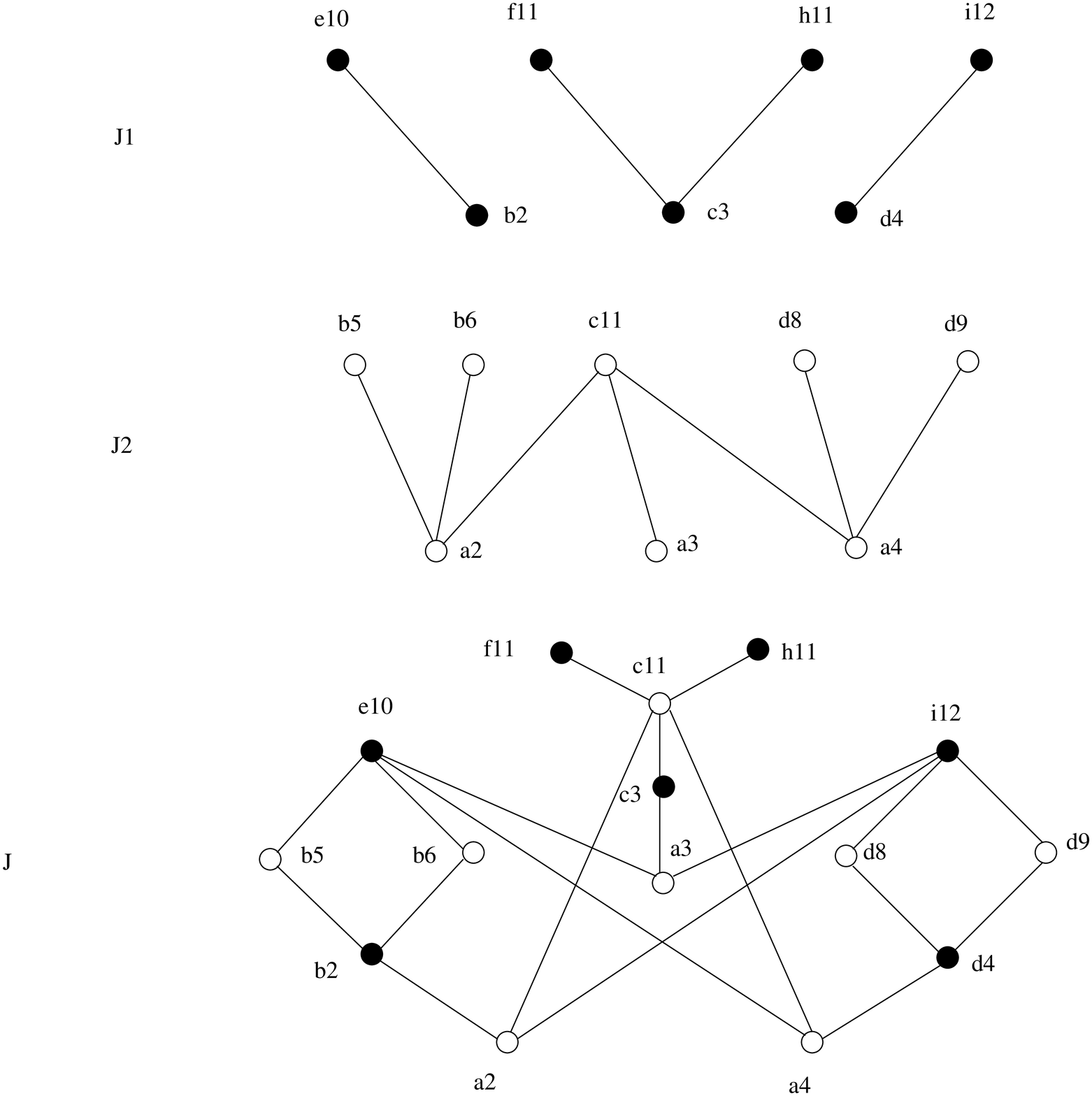}
\caption{}
\label{fig4}
\end{center}
\end{figure}

Our task to generate $L$ as the closure system $\mathscr{C}=\{J(x)\;|\;x\in L\}$ of all $\Lambda$-closed order ideals $X$ of $(J,\le)$ amounts to determine those $X=Y\cup Z\; (Y\in J_1[\Lambda _1], Z\in J_2[\Lambda _2])$ that happen to be order ideals of $(J,\le)$. For $X=Y\cup Z$ to be an order ideal of $(J,\le)$ it is necessary  (but not sufficient) that $Y$ and $Z$ be order ideals of $(J_1,\le)$ and $(J_2,\le)$ respectively. In particular, since $d4\le i12$ in $(J_1,\le)$ and since the fourth row from below in table \ref{tab1} has $0$ at position $d4$, the $i12$-component $\varepsilon _2$ can safely be switched to $0$. Accordingly the other $\varepsilon _2$ turns to be $2$. (If there had been $1$ at position $i12$, then the row must have been deleted.) This yields $r_1$ below.  Together with three similarly obtained rows we get a certain subset $J^*_1[\Lambda _1]$ of $J_1[\Lambda _1]:$
\newpage
\begin{table}[!h]
\begin{center}
\begin{tabular}{l|l|l|l|l|l|l|l|}
  & $b2$ & $c3$ & $d4$ & $f11$ & $h11$ & $i12$ & $e10$\\
 \hline
 $r_1 =$& $ \varepsilon$ & $\varepsilon $ & $ 0$ & $2 $ & $ 0$ & {\bf 0} & {\bf 2}\\
 \hline
 $r_2 =$& $0 $ & $ 1$ & $0 $ & $0 $ & $ 1$ & $ 0$ & $ 0$ \\
 \hline
 $r_3 =$& $ \delta$ & $\delta $ & $1 $ & $0 $ & $0 $ & $\varepsilon  $ & $\varepsilon $ \\
 \hline
 $r_4 =$& $1 $ & $ 1$ & $ 1$ & $ 1$ & $1 $ & $\delta $ & $\delta $ \\
 \hline
\end{tabular}
\caption{$J^*_1[\Lambda _1]$}
\label{tab2}
\end{center}
\end{table}

Computing $J_2[\Lambda _2]$ along the same lines, and again filling in the "immediate" $0$'s and $1$'s forced by the poset $(J_2,\le)$, yields this subset $J^*_2[\Lambda _2]$ of $J_2[\Lambda _2]$
(where $\delta \delta \delta$ is "$000$ or $111$"):

\begin{table}[!h]
\begin{center}
\begin{tabular}{l|l|l|l|l|l|l|l|l|}
  & $a2$ & $a3$ & $a4$ & $b6$ & $b5$ & $c11$ & $d9$ & $d8$\\
 \hline
 $s_1 =$& $ \varepsilon$ & $\varepsilon $ & $ 0$ & $2 $ & $ 0$ & $2$ & $0$ & $0$\\
 \hline
 $s_2 =$& $1$ & $ 0$ & $0 $ & $0 $ & $ 1$ & $ 0$ & $ 0$& $0$ \\
 \hline
 $s_3 =$& $ 0$ & $0$ & $1 $ & $0 $ & $0 $ & $0  $ & $0 $ & $2$\\
 \hline
 $s_4 =$& $0$ & $ 0$ & $ 1$ & $ 0$ & $0 $ & $0$ & $1 $ & $0$ \\
 \hline
 $s_5 =$& $1$ & $ 1$ & $1 $ & $0 $ & $ 0$ & $ 2$ & $ 0$ & $0$ \\
 \hline
 $s_6 =$& $ 1$ & $1 $ & $1 $ & $0 $ & $0 $ & $0 $ & $1 $ & $1$ \\
 \hline
 $s_7 =$& $1 $ & $ 1$ & $ 1$ & $ 1$ & $1 $ & $\delta $ & $\delta $ & $\delta$\\
 \hline
\end{tabular}
\caption{$J^*_2[\Lambda _2]$}
\label{tab3}
\end{center}
\end{table}
Each member of $\mathscr{C} \simeq L$, i.e. each $\Lambda$-closed order ideal $X$ of $(J,\le)$, is of the form $X=Y\cup Z$ for some $i\in \{1,\cdots, 4\}$, $Y\in r_i$ and for some $j\in\{1,\cdots,7\}$, $Z\in s_j$. In order to get these $X$'s we first discard the "concatenated" rows $r_is_j$ which do not contain {\it any} order ideal of $(J,\le)$. For instance $r_1s_4$ is of that kind: Because\footnote{We can read $d4 < d9$ from Figure \ref{fig4}, the algorithm "knows" it from the given implication $\{d9\}\ra J(d9)$.} $d4 < d9$, no  order ideal $X$ has $1$ at $d9$ but $0$ at $d4$. We say that $1$ and $0$ {\it clash}. As another example, suppose $X$ was an order ideal contained in $r_2s_1$. From $c3, h11 \in X$ (see $r_2$) and $c3>a3$ and $h11>a2$ follows $a3,a2\in X$. But this cannot be since the $\varepsilon \varepsilon $ in $s_1$ forces $\left|X\cap \{a3,a2\}\right|\le 1$. Hence also $r_2s_1$ contains no order ideals. 

In this way one finds that at most
\begin{eqnarray}\label{eq9}
r_1s_1,\; r_1s_2,\; r_1s_3,\; r_1s_5,\; r_1s_7,\; r_2s_5,\;  r_3s_3,\; r_3s_4,\; r_3s_5,\; r_3s_6,\; r_3s_7,\; r_4s_5,\; r_4s_7
\end{eqnarray}
contain order ideals. In order to filter them from each of these 13 concatenated rows we impose all implications $\{p\}\ra J(p)$ with the $(a,B)$-algorithm from section $2$ in Part I. Of course, instead of $J(p)$ it suffices to take the smaller set of lower covers of $p$. For singleton premises (as in $\{a\}\ra B$) the implication $n$-algorithm can be streamlined to the $(a,B)$-algorithm discussed in \cite{mw1}.

Consider e.g.
 $$r_1s_5=(\varepsilon,\varepsilon, 0,2,0,0,2,1,1,1,0,0,2,0,0)$$
which we shall work from left to right. The implication $\{b2\}\ra \{a2\}$ holds already since each $X\in r_1s_5$ has $a2\in X$. Similarly for $\{c3\} \ra \{a3\}$ and $\{d4\}\ra \{a4\}$. As to $\{f11\}\ra \{c11\}$, the corresponding components are both $2$ and hence can be turned to $a,b$ respectively. The next not yet holding implication is $\{e10\}\ra \{b5,b6,a3\}$. Since say $b6\notin X$ for all $X\in r_1s_5$, we turn $2$ to $0$ at position $e10$. So far we have 
$$\rho =(\varepsilon, \varepsilon, 0,a,0,0,0,1,1,1,0,0,b,0,0).$$
 In order to impose the next not yet holding implication, i.e. $\{c11\}\ra \{a2,a4,c3\}$, we need to split $\rho$ as follows:
\begin{eqnarray*}
\rho _1&= & (\varepsilon, \varepsilon, 0,0,0,0,0,1,1,1,0,0,{\bf 0},0,0)\\
\rho _2&= & (0, 1, 0,2,0,0,0,1,1,1,0,0,{\bf 1},0,0)
\end{eqnarray*}
Notice that {\bf 0} in $\rho _1$ turns {\it $a$} to $0$ at position $f11$. Further, {\bf 1} in $\rho _2$ turns {\it $a$} to $2$ at position $f11$, and $\varepsilon$ to $1$ at position $c3$. Hence the other $\varepsilon$ becomes $0$. The other concatenated rows in (\ref{eq9}) are treated similarly, and the result is this:

\begin{longtable}{l|l|l|l|l|l|l|l|l|l|l|l|l|l|l|l|}
  & $b2$&$c3$& $d4$& $f11$ & $h11$ & $i12$ & $e10$& $a2$& $a3$& $a4$ &$b6$ &$b5$& $c11$& $d9$& $d8$\\
  \hline
    & $0$&$0$& $0$& $0$ & $0$ & $0$ & $0$& $\varepsilon$& $\varepsilon$& $0$ &$0$ &$0$& $0$& $0$& $0$\\
    & $0$&$1$& $0$& $0$ & $0$ & $0$ & $0$& $0$& $1$& $0$ &$0$ &$0$& $0$& $0$& $0$\\
    & $1$&$0$& $0$& $0$ & $0$ & $0$ & $0$& $1$& $0$& $0$ &$2$ &$0$& $0$& $0$& $0$\\
     & $1$&$0$& $0$& $0$ & $0$ & $0$ & $0$& $1$& $0$& $0$ &$0$ &$1$& $0$& $0$& $0$\\
      & $0$&$0$& $0$& $0$ & $0$ & $0$ & $0$& $0$& $0$& $1$ &$0$ &$0$& $0$& $0$& $0$\\
 $\rho _1\,=$  & $\varepsilon$&$\varepsilon$& $0$& $0$ & $0$ & $0$ & $0$& $1$& $1$& $1$ &$0$ &$0$& $0$& $0$& $0$\\
 $\rho _2 \,=$ & $0$&$1$& $0$& $2$ & $0$ & $0$ & $0$& $1$& $1$& $1$ &$0$ &$0$& $1$& $0$& $0$\\
     & $1$&$0$& $0$& $0$ & $0$ & $0$ & $2$& $1$& $1$& $1$ &$1$ &$1$& $0$& $0$& $0$\\
     & $0$&$1$& $0$& $0$ & $1$ & $0$ & $0$& $1$& $1$& $1$ &$0$ &$0$& $1$& $0$& $0$\\
     & $0$&$0$& $1$& $0$ & $0$ & $0$ & $0$& $0$& $0$& $1$ &$0$ &$0$& $0$& $0$& $2$\\
     & $0$&$0$& $1$& $0$ & $0$ & $0$ & $0$& $0$& $0$& $1$ &$0$ &$0$& $0$& $1$& $0$\\
     & $\delta$&$\delta$& $1$& $0$ & $0$ & $0$ & $0$& $1$& $1$& $1$ &$0$ &$0$& $0$& $0$& $0$\\
     & $1$&$1$& $1$& $0$ & $0$ & $0$ & $0$& $1$& $1$& $1$ &$0$ &$0$& $1$& $0$& $0$\\
     & $\delta$&$\delta$& $1$& $0$ & $0$ & $2$ & $0$& $1$& $1$& $1$ &$0$ &$0$& $0$& $1$& $1$\\
     & $1$&$1$& $1$& $0$ & $0$ & $0$ & $2$& $1$& $1$& $1$ &$1$ &$1$& $\delta$& $\delta$& $\delta$\\
     & $1$&$1$& $1$& $0$ & $0$ & $1$ & $0$& $1$& $1$& $1$ &$1$ &$1$& $1$& $1$& $1$\\
     & $1$&$1$& $1$& $1$ & $1$ & $0$ & $0$& $1$& $1$& $1$ &$0$ &$0$& $1$& $0$& $0$\\
     & $1$&$1$& $1$& $1$ & $1$ & $0$ & $0$& $1$& $1$& $1$ &$1$ &$1$& $1$& $1$& $1$\\
     & $1$&$1$& $1$& $1$ & $1$ & $1$ & $1$& $1$& $1$& $1$ &$1$ &$1$& $1$& $1$& $1$\\
     \hline
     \caption{}
     \label{tab4}
\end{longtable}
Exactly the 34 $\Lambda$-closed order ideals $X=J(x)\;(x\in L)$ are encoded in this table. For instance, letting $\varepsilon \varepsilon =01$ in row $\rho _1$ we get $X=\{c3, a2, a3, a4\}$ which is $J(c7)$ (see Fig.2).

Let us recap the described procedure. Steps (d),(e) and (f) convey an extension of the method that pays off for large $t$.

\textit{ \large{Summary: Calculating a modular subdirect product from the connection maps}}. 

\begin{itemize}
\item[(a)] For each factor lattice $L_i$ find a base of lines $\Lambda'_i\; (1\le i\le t)$.
\item[(b)] Calculate the connected components $(J_i,\Lambda_i)$ where $J_i:=\sigma _i\left(J(L_i)\right)$ and $\Lambda _i:=\sigma _i\left(\Lambda'_i\right)$. Here $\sigma _i: L_i \ra \ds \prod^t_{j=1}L_j$ is calculated as $\sigma _i(x):=(\alpha _{j,i}(x):\; 1\le j\le t)$.
\item[(c)] Using the described $\delta, \varepsilon, l$-algorithm (more details in \cite{mw1}) compute a context of each $J_i[\Lambda _i]$, i.e. compute a compact representation for the family of all $\Lambda _i$-closed subsets of $J_i \;(1\le i\le t)$. Advantageous, but not strictly necessary are certain subfamilies $J^*_i[\Lambda _i]\subseteq J[\Lambda _i]$  (as in table \ref{tab2}, table \ref{tab3}) because they have shorter contexts, which reduces the size of a same graph $G$ in the next step.
\item[(d)] Let $J:=\ds \bigcup _{1\le i\le t}J_i$. As subset of the known lattice $\ds \prod_{1\le i\le t}L_i$ the set $J$ becomes partially ordered. Consider the graph $G$ whose vertices are the rows occuring in the contexts  $J^*_i[\Lambda _i]\;(1\le i\le t)$.  Let $k_i$ be the number of rows of $J_i^*[\Lambda_i]$. By definition these $t$ contexts constitute disjoint $k_i$-cliques of $G$. Moreover, two rows from distinct 
 cliques are declared adjacent if they contain components $1$ and $0$ respectively that \textit{clash} (with respect to the partial ordering $(J,\le)$, as seen in the $t=2$ example).
\item[(e)] Calculate {\it all} $t$-element anticliques of $G$, for instance with the algorithm of \cite{mw3}.
\item[(f)] The rows $\rho$ concatenated from the $k_1k_2\cdots k_t$ transversals  $\{\rho _i\;|\; 1\le i\le t\}$ of the $t$ contexts comprise precisely the $\Lambda$-closed subsets of $J$. Such  a row $\rho$ is {\it good} in the sense of containing at least one order ideal $X$ of $(J,\le)$ if and only if $\{\rho _i\;|\; 1\le i\le t\}$ is an anticlique of $G$; and the latter have been computed in (e). Using the $(a,B)$-algorithm to impose all implications $\{p\}\ra J(p)\; (p\in J)$  on a good row $\rho$ filters the order ideals $X$ from it.
\end{itemize}

\section{Application to lattices freely generated by posets within  f.g. varieties}\label{section6}

Recall that a variety ${\cal V}$ of lattices is {\it finitely generated} if it is generated by a single finite lattices. Equivalently, and more to the point for us, ${\cal V}$ has up to isomorphism only finitely many subdirectly irreducibles $S_1,\cdots, S_r$, and they are all finite. Thus every $L \in {\cal V}$ is a subdirect product of lattices $L_i\;(1\le i\le t)$ where each $L_i$ is isomorphic to some $S_k$. Possibly $L_i \simeq L_j\simeq S_k$ for $i\ne j$.

Let $(P,\le)$ be a finite poset. We wish to compute the lattice $F{\cal V}(P,\le)$ freely generated by $(P,\le)$ within ${\cal V}$ as defined in Part I. 
If we  knew the precise structure of the connection maps $\alpha_{ij}: L_j \ra L_i$, then we could construct the subdirect product $F\mathcal{V}(P,\le)$ as in section \ref{section2} and \ref{section3}! 

This works out as follows. Restricting the projections $F\mathcal{V}(P,\le)\ra L_i$ to $P$ yields a {\it $P$-labelling} of $L_i$, i.e. a monotone map $\lambda$ from $P$ onto a generating set of $L_i$. Conversely, by the universal mapping property, {\it each} $P$-labelling arises in this way. We will use the following poset as a standard poset as in Part I:

\begin{figure}[!h]
\begin{center}
\psfrag{P}{$(P,\le)\;\;=$}
\includegraphics[scale=.4]{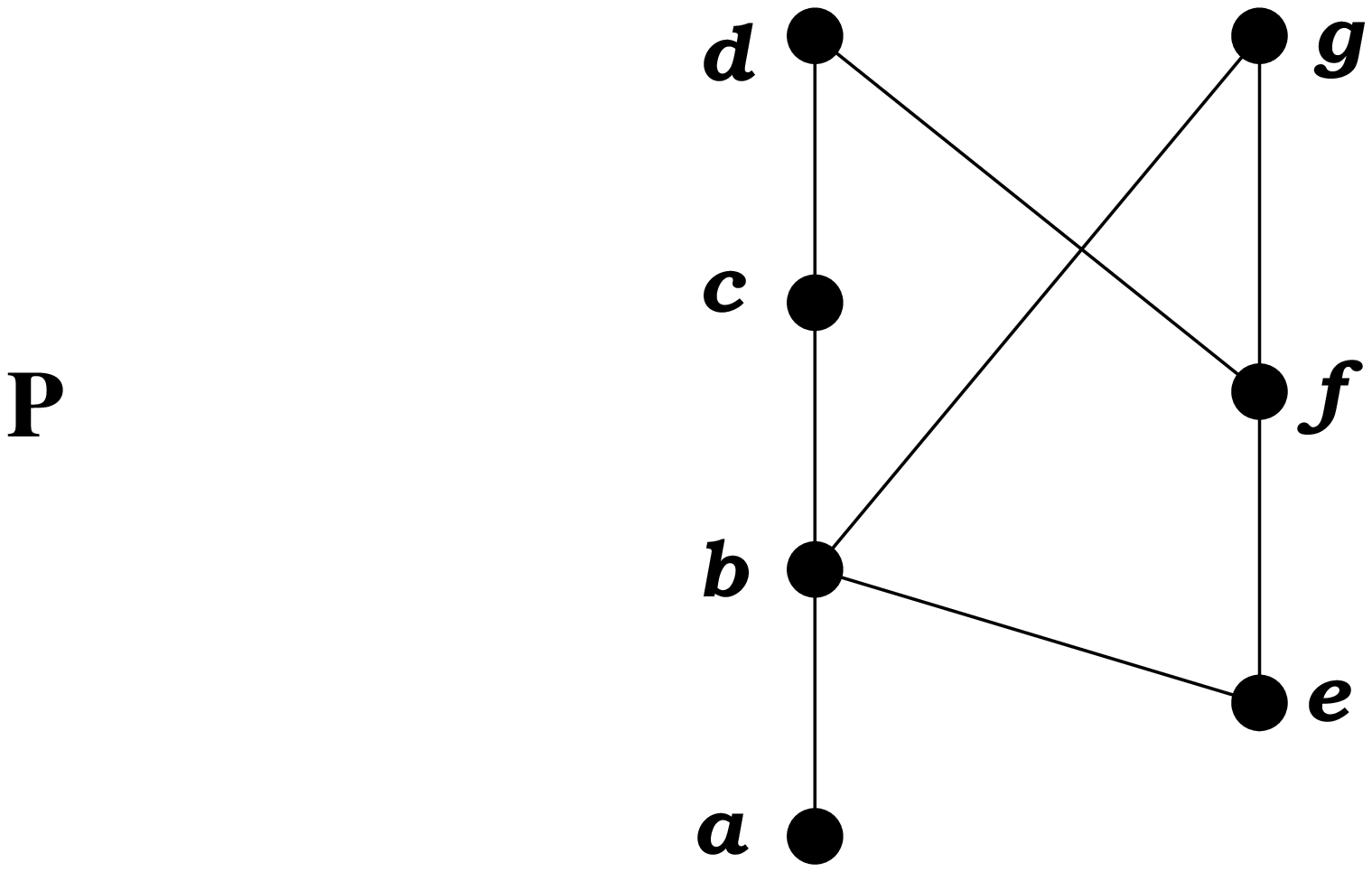}
\caption{}
\label{standardposet}
\end{center}
\end{figure}

\begin{ex} \textrm{The finitely generated  variety} $\mathcal{V} = \mathcal{V}( _5N_5)$ \textrm{has} $r = 3$ \textrm{subdirectly irreducibles} $S_1, S_2, S_3$ \textrm{which (renamed) are these:}
\end{ex}

\begin{figure}[!h]
\begin{center}
 \includegraphics[scale=.4]{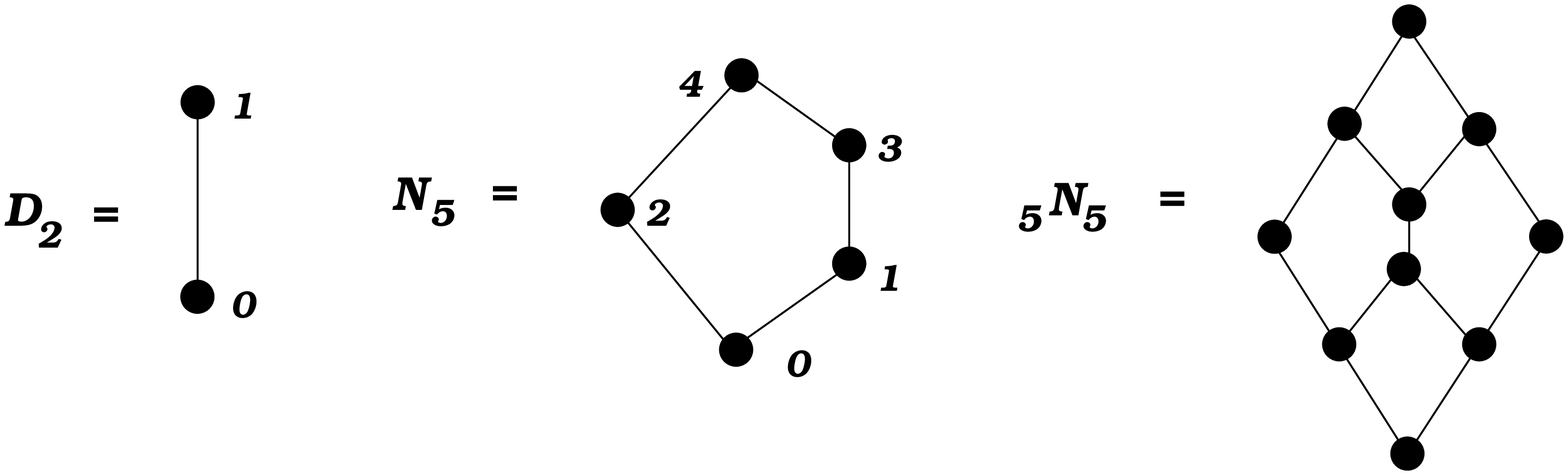}
\caption{}
\end{center}
\end{figure}

If $P$ is the poset of Figure \ref{standardposet}, what are the $P$-labellings of these lattices? The twelve $P$-labellings $\lambda$ of $D_2$ are monotone maps and hence the sets $\lambda ^{-1}(1)$ yield  the nonempty filters of $(P, \leq)$. We encountered the twelve $P$-labellings of $D_2$ already in Part I but we computed the free distributive lattice $FD(P,\le)$ by other means.

These are the seven $P$-labellings of $N_5$ (for readability $0,1,2,3,4$ are written only once):
\newpage
\begin{figure}[!h]
\begin{center}
 \includegraphics[scale=.425]{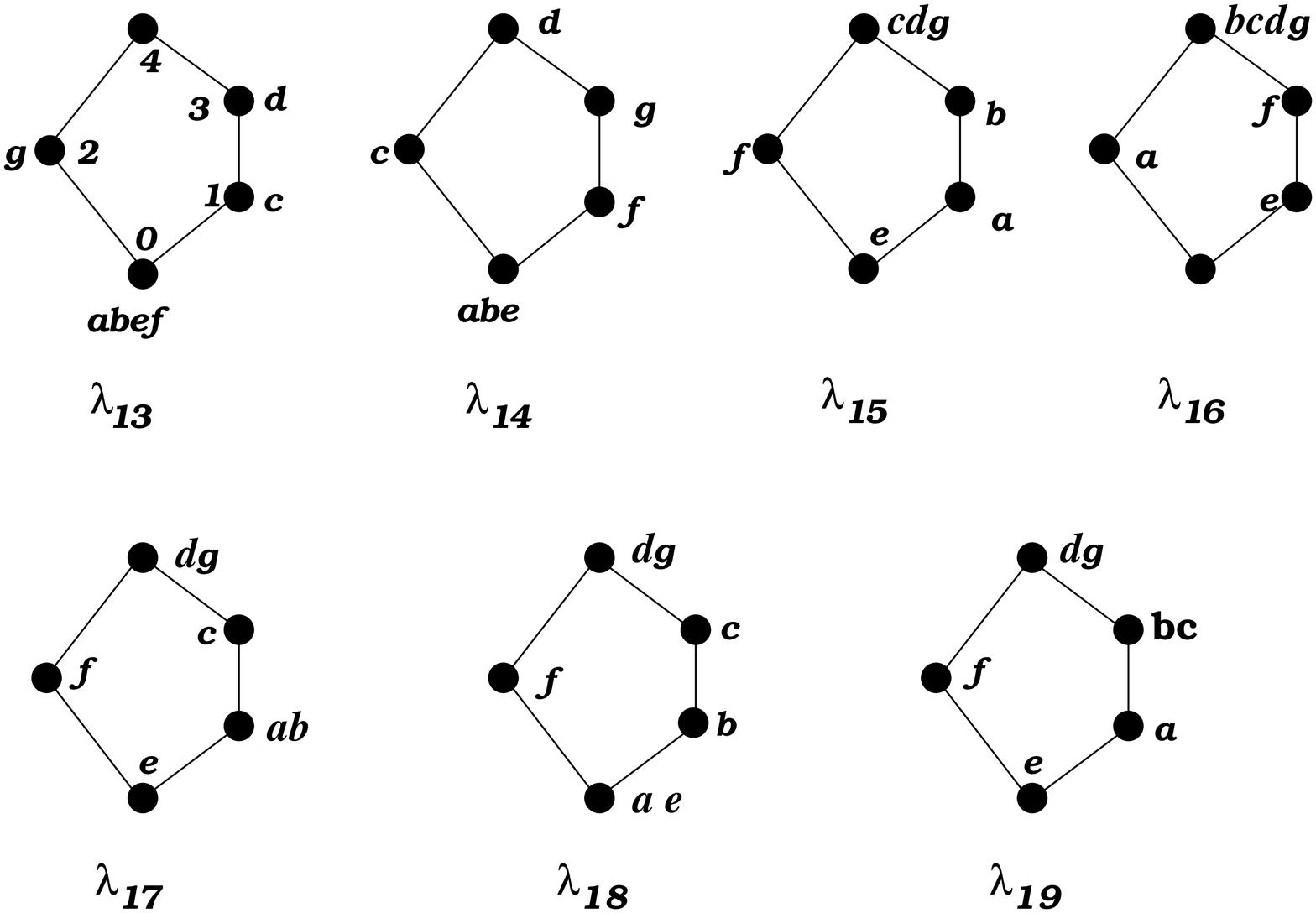}
\caption{}
\label{fig8}
\end{center}
\end{figure}

Here are two $P$-labellings of $_5N_5$:

\begin{figure}[!h]
\begin{center}
\includegraphics[scale=.45]{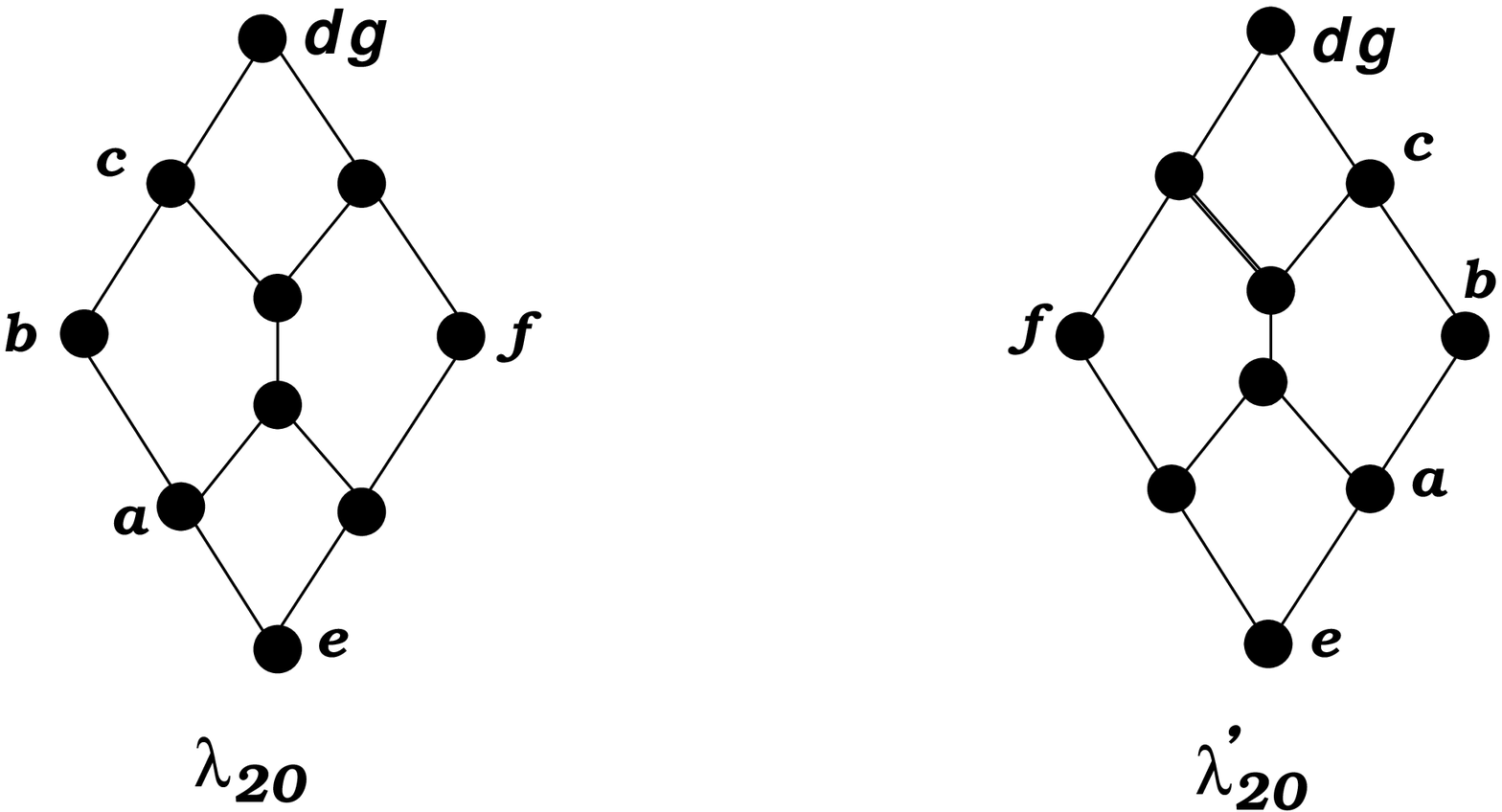}
\caption{}
\label{fig9}
\end{center}
\end{figure}

They are {\it equivalent} in the sense that   $\lambda'_{20} = \alpha \circ \lambda_{20}$ for some automorphism $\alpha$ of $_5N_5$. One can show that up to automorphism $\lambda_{20}$ is the only $P$-labelling of $_5N_5$. Therefore $F \mathcal{V} (P,\le)$ is a certain subdirect product of  $20$ lattices $L_i$, twelve of which are isomorphic to $D_2$, seven to $N_5$, and one to $_5N_5$.

It turns out, crucially, that the connecting $\vee$-morphism $\alpha_{i,j}: L_j \ra L_i$ is the {\it biggest} $\vee$-homomorphism that maps labels below corresponding labels [3, Satz 3.6].  For instance,  between $L_{14}$ and $L_{17}$ we have
\newpage 

\begin{figure}[!h]
\begin{center}
 \includegraphics[scale=.45]{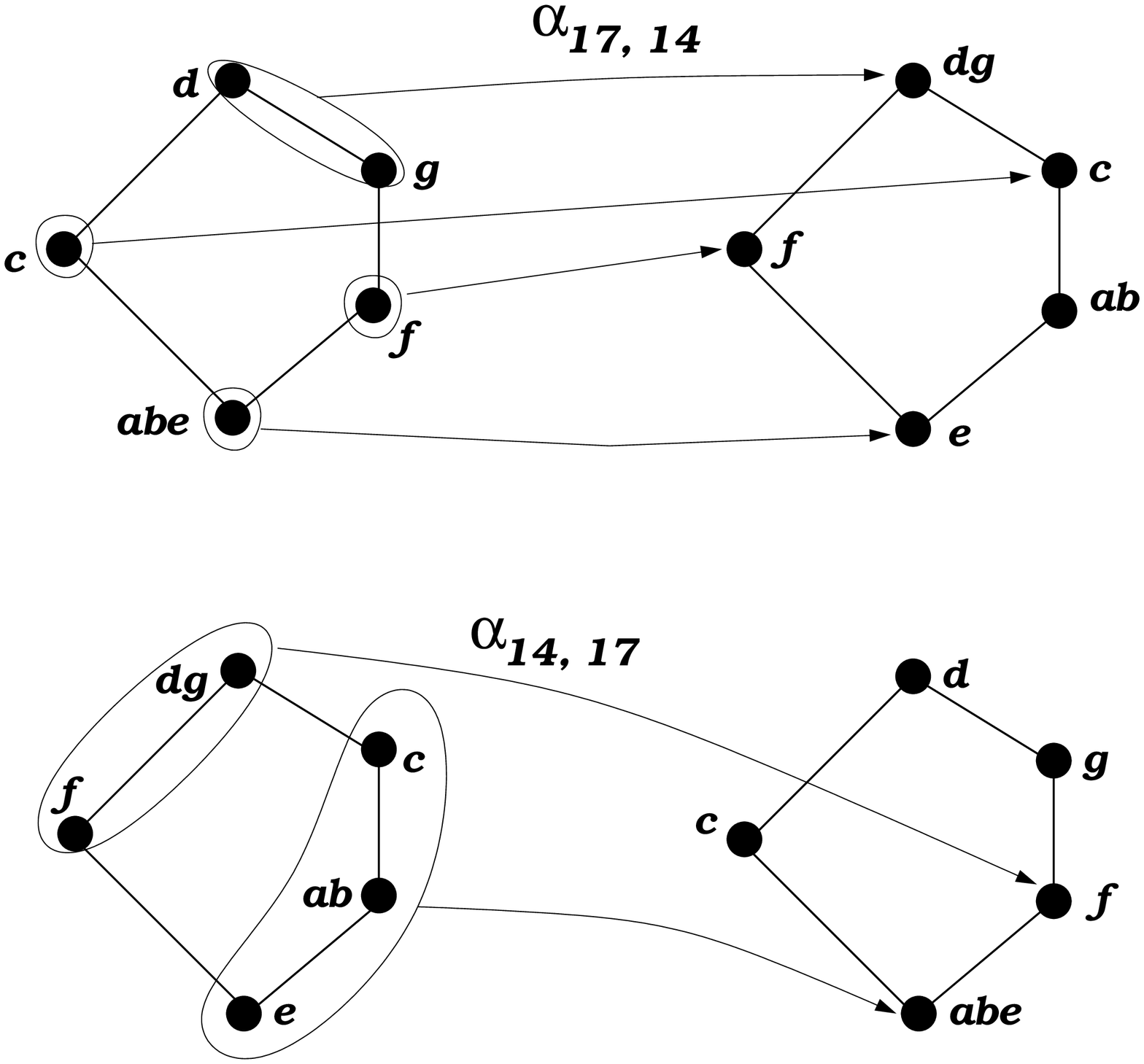}
\caption{}
\label{fig10}
\end{center}
\end{figure}

Notice that the $\vee$-homomorphism $\alpha = \alpha_{17,14}$ does not respect meets,  $\alpha (c \wedge g) \neq \alpha (c) \wedge \alpha (g)$, but $\alpha_{14, 17}$ happens to be a lattice homomorphism. 

Let $F\mathcal{V}_1(P,\le), F\mathcal{V}_2(P,\le)$, and $F\mathcal{V}_3(P,\le)$ be the factor lattices of $F\mathcal{V}(P,\le)$ obtained by taking the subdirect products of $L_1, \cdots, L_{12}$ respectively $L_{13}, \cdots, L_{19}$, respectively $L_{20}$. These {\it homogeneous components} of $F\mathcal{V}(P,\le)$ are $F\mathcal{V}_3(P,\le)\simeq {}_5N_5$ and $F\mathcal{V}_1(P,\le)\simeq FD(P,\le)$; as well as  $F\mathcal{V}_2(P,\le)$ which is depicted below
 together with its previously mentioned factor lattices $L_{14}$ and $L_{17}$:

\newpage
\begin{figure}[!h]
\begin{center}
 \includegraphics[scale=.45]{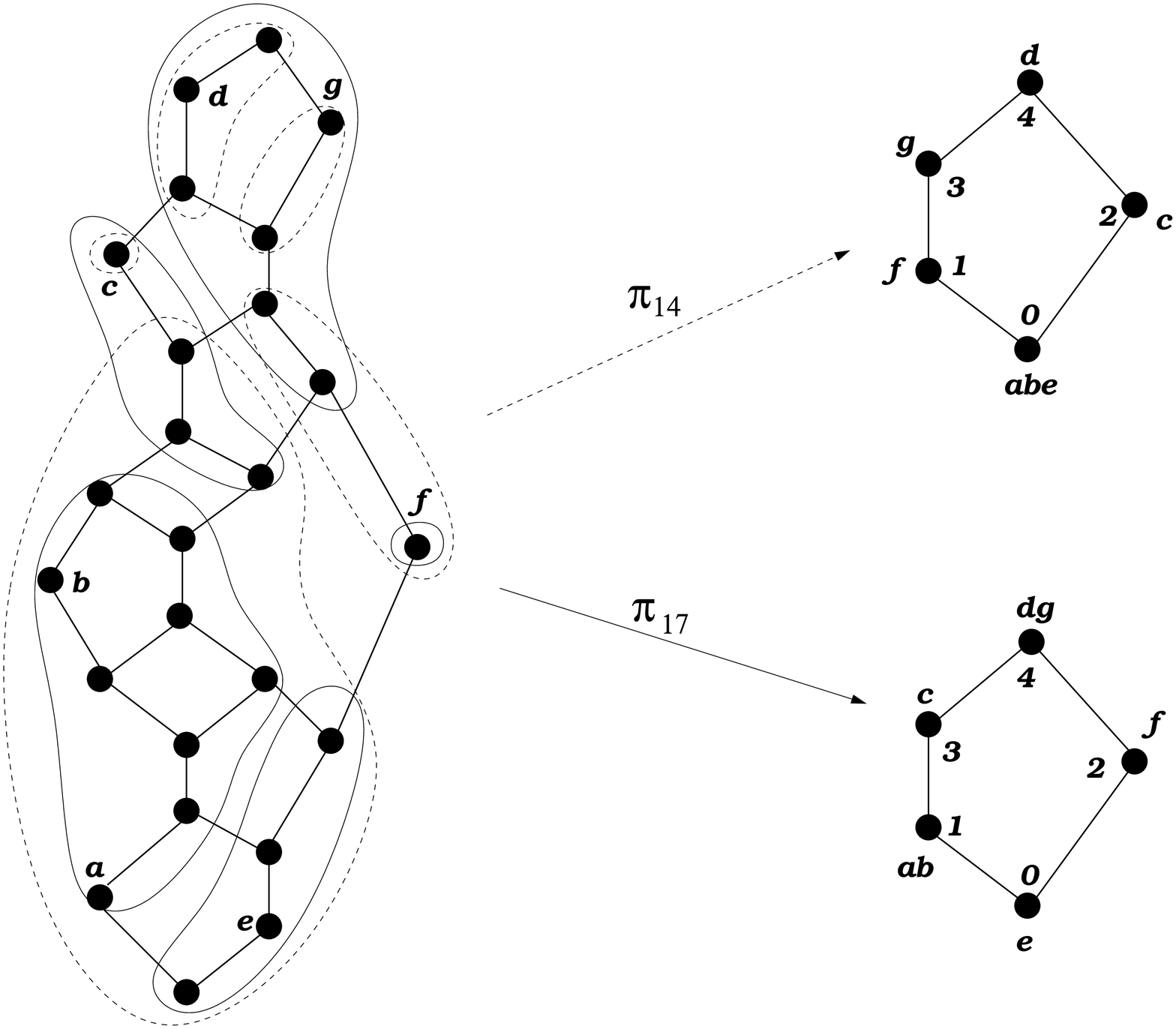}
\caption{}
\label{fig11}
\end{center}
\end{figure}

Observe that say $\pi_{14} \circ \alpha_{17}$ coincides indeed with the initially given $\alpha_{14, 17}$. For instance,
$$(\pi_{14} \circ \alpha_{17}) (2) \;\; = \;\; \pi_{14} (f) \;\;=\;\;1\;\;=\;\;\alpha_{14, 17} (2).$$

From the diagram above it is evident that the join irreducible $f \in F\mathcal{V}_2(P,\le)$ belongs to both \footnote{The reader may check that also $f \in \alpha_{18} (2) \cap \alpha_{19}(2)$.} $\alpha_{14}(1)$ and $\alpha_{17}(2)$, which illustrates that the union in (\ref{eq1}) needs not be disjoint for non-modular lattices $L$.

Having the connection maps $\alpha_{i,j} \;(1\le i,j\le 20)$ at hand, the lattice $F:=F\mathcal{V}(P,\le)$ can be computed as discussed in section \ref{section2}. Specifically, recall that the maps $\alpha_{i,j}$ yield the maps $\sigma'_i$ from $L_i$ to $L_1\times \cdots \times L_{20}$ according to (\ref{eq3}). The $\sigma'_i$ in turn yield the subsets $J(F)$ and $G(F)$ of $L_1\times \cdots \times L_{20}$ according to (\ref{eq1}) and (\ref{eq5}). Now $F$ can be computed by running the $(A,B)$-algorithm on all implications $A\ra J(\bigvee A)$ where $A\subseteq J(F)$ is such that $\bigvee A\in G(F)$.

The above method works to compute $F\cal{V}(P,\le)$ for {\it any} finitely generated variety $\cal{V}$ of lattices. Our particular choice of $\cal{V}=\cal{V}$$({}_5N_5)$ was motivated by some extra feature of this variety. Namely, recall that $FL(P,\le)$ is the lattice freely generated by $(P,\le)$ within the (not f.g.) variety of {\it all} lattices. It turns out that  when $FL(P,\le)$ happens to be finite, it coincides with $F\mathcal{V}(P,\le)$. Finiteness takes place \cite{Wi} if and only if $P$ has no subposet isomorphic to ${\bf 1 + 1+ 1}$ or ${\bf 2 + 2}$ or $ {\bf 1 + 4}$.  Here, say ${\bf 2+2}$ denotes the disjoint union of two $2$-element chains. 

If partial semilattices $(P, \vee')$ rather than mere posets $(P, \leq)$ are at stake, everything "should" stay the same, except that there are usually less $(P, \vee')$-{\it labellings} $\lambda_j$. The latter by definition are not just monotone but also respect the declared suprema, i.e. $\lambda_j(a\vee' b)=\lambda_j(a)\vee \lambda_j(b)$. Detailed proofs are still pending.

\subsection{A symmetry exploiting variation}

 Consider the natural epimorphism
\[
	f: F{\cal V}(P,\le)\ra FD(P,\le)
\]
The idea is to calculate $F{\cal V}(P,\le)$ as the disjoint union of the interval sublattices $f^{-1}(x)$ with $x$ ranging over $x\in FD(P,\le)$. Before going into further details, notice that this approach is appealing when $(P,\le)$ has a large automorphism group $\cal{S}$ and hence decays into few and large $\cal{S}$-orbits $\Omega _i.$ This is because $f^{-1}(x)\simeq f^{-1}(y)$ for all $x, y\in \Omega _i$, and thus only one $f^{-1}(x)$ per orbit needs to be computed.

As to the computation of $K:= f^{-1}(x)$, let $F{\cal V}(P,\le)$ be a subdirect product of the subdirectly irreducible lattices $L_i$ and let $\phi _i: \ds \prod_{1\le j\le t}\!\!\!L_j \ra L_i$ be the canonical projections, and $\sigma _i: L_i\ra \ds \prod_{1\le j\le t}\!\!\!L_j$ the corresponding (known) smallest pre-image maps. Setting $K_i =\phi _i(K)\subseteq L_i$ it is clear that $K$ is the subdirect product of the lattices $K_i \;(1\le i\le t)$.  When the sublattices $K_i$ of $L_i$ are known, one  can calculate the scaffolding
\begin{eqnarray*}
G\left(K\right) &=& \bigcup_{1\le i\le t}\sigma _i \left(K_i\setminus \{0\}\right)
\end{eqnarray*} 
Here the maps $\sigma _i$ are still the {\it same} as for $L_i$. Using the $(A,B)$-algorithm or variations thereof  the lattice $K$ can then be computed as the  closure system of all $\bigvee$-ideals of the partial semilattice $\left(G(K),\bigvee \right)$.

But how is $K_i$ computed? Each $x$  in $FD(P,\le)$ can be written, in many ways, as a lattice polynomial of elements of $P$. Considered within $F{\cal V}(P,\le)$ some of these lattice polynomials may yield distinct elements. Let $DNF(x)$ be the unique disjunctive normal form of $x$, and identify $DNF(x)$ with the corresponding element in $F{\cal V}(P,\le)$. Similarly define $CNF(x)$ in terms of the conjunctive normal form. To fix ideas, say $a,b,c,d \in P \subseteq FD(P,\le)$ and the corresponding elements in $F{\cal V}(P,\le)$ are $a',b',c',d'$. If $x=\left((a\vee b)\wedge c\right)\vee d$, then 
\begin{eqnarray*}
DNF(x)& =&(a'\wedge c')\vee (b'\wedge c')\vee d' \;\,\le \;\,\left((a'\vee b')\wedge c'\right)\vee d'\;\,\le\;\, (a'\vee b'\vee d')\wedge (c'\vee d')\\  &= & CNF(x).
\end{eqnarray*}
Provided that $(P,\le)$ is unordered (an antichain), it is shown in [1, Thm.3.3] that for all $x\in FD(P,\le)$ one has
\begin{eqnarray*}
K &=& f^{-1}(x) \;\;\;=\;\;\; \left[DNF(x), CNF(x)\right]
\end{eqnarray*}
Thus, for $x$ as above we get $K_i =\left[\phi _i\left(DNF(x)\right), \phi _i\left(CNF(x)\right)\right]$ with e.g. 
\begin{eqnarray*}
\phi _i\left(DNF(x)\right) &=& \left(\phi _i(a')\vee \phi _i(b')\vee \phi _i(d')\right)\wedge \left(\phi _i(c')\vee \phi _i(d')\right).
\end{eqnarray*}
This is readily evaluated because $\phi _i(a'),\cdots, \phi _i(d')$ are just some of the known labels of $L_i$. When $(P,\le)$ is not an antichain, Thm.3.3 in \cite{berman} needs to be adapted. Probably this is easy. 

\section{The smallest modular non-distributive variety} \label{section7}

The smallest modular nondistributive variety  is $\cal{V}=\cal{V}$$(M_3)$, and it has $D_2$ and $M_3$ as subdirectly irreducibles. Since $\cal{V}$ is finitely generated, the computation of $FM_3(P,\le):=F\cal{V}(P,\le)$ works according to section $6$. On the other hand, $FM_3(P,\le)$ is a modular lattice, and so the specialities of section \ref{section5} apply.

Specifically, in step $(a)$ at the end of section \ref{section5} each $\Lambda'_i$ merely consists of \textit{one} $3$-element line. Step $(b)$ involves the calculation of all $P$-labellings of $D_2$ and $M_3$, as well as the biggest $\vee$-morphisms $\alpha_{j,i}$ that map labels below corresponding labels.  The explicit programming of all of that was done with Mathematica. Steps $(c)$ to $(f)$ were  condensed considerably, but  for finitely generated modular varieties with more or bigger subdirectly irreducibles these steps would presumably pay off.

Our variety $\cal{V}$$(M_3)$ enjoys an extra property akin to the variety $\cal{V}$$(_5N_5)$ in section $6$. That is, whenever the free {\it modular} lattice $FM(P,\le)$ generated by the finite poset $P$ happens to be finite, then $FM(P,\le)$ coincides with $FM_3(P,\le)$. Finiteness takes place if and only if $P$ has no subposet isomorphic to {\bf 1+1+1+1} of {\bf 1+2+2}. All of this is due to Wille 1973. We mention that an English version, and also a more explicit one with helpful drawings, of Wille's German proof, features in \cite{yves}. 

For all $1+2+5+16+63 =87$ posets $P$ with $|P|\le 5$ the lattices $FD(P,\le)$ have been drawn (some in compressed form) in \cite{bartenschlager}. In \cite{yves} their cardinalities were recalculated and confirmed. In fact the cardinalities of $FD(P,\le)$ {\it and} $FM_3(P,\le)$ are calculated in \cite{yves} for almost all posets $P$ with $|P|\le 6$, and the numbers $s$ and $s+t$ of subdirectly irreducible factors of $FD(P,\le)$ respectively $FM_3(P,\le)$ are listed. Observe that $FM_3(P,\le)$ has length $s+2t$. As to "almost", the list lacks $14$ out of $318$ six element posets which due to their high symmetry blew up  $FM_3(P,\le)$ too much. However, chances are good that implementing the symmetry exploiting ideas of $6.1$ would finish the job.  We mention that the cardinalities
\[
28, \;\;\;138, \;\;\;629, \;\;\;2784
\]
of {\bf 1+1+1} (the Dedekind lattice), {\bf 1+1+2} (the so called Takeuchi lattice), {\bf 1+1+3}, and {\bf 1+1+4} match the explicit formula for $|FM_3(1+1+n)|$ found in \cite{peter}.

Here is the data for the first few $6$-element posets. 
\begin{longtable}{|c|c|c|c|c|c|c|c|}
\hline
\includegraphics[scale=.25]{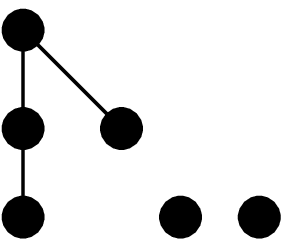} & 1326 &\textbf{{\small 296198143}}& 26+45&  \includegraphics[scale=.25]{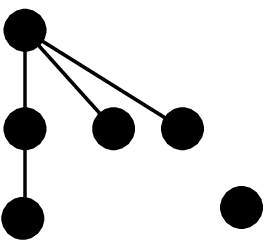} & 936 &\textbf{{\small 160224000}}& 24+39\\
\hline
 \includegraphics[scale=.25]{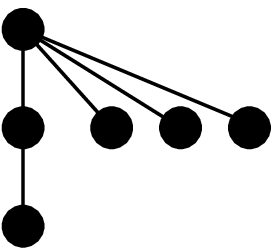} & 886 &\textbf{{\small 160228750}}& 23+39&  \includegraphics[scale=.25]{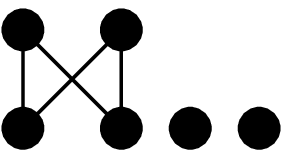} & 1058 &\textbf{{\small 6306868}}& 26+37\\
\hline
 \includegraphics[scale=.25]{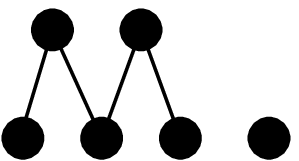} & 670 &\textbf{{\small 434366}}& 24+26&  \includegraphics[scale=.25]{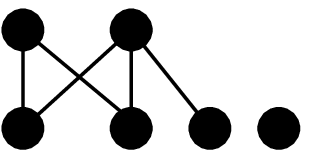} & 407 &\textbf{{\small 68915}}& 22+20\\
\hline
 \includegraphics[scale=.25]{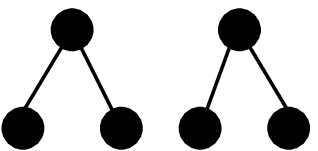} & 590 &\textbf{{\small 2472286}}& 23+22&  \includegraphics[scale=.25]{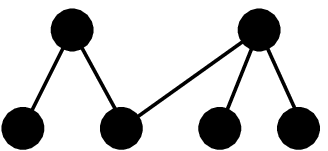} & 354 &\textbf{{\small 64461}}& 21+18\\
\hline
 \includegraphics[scale=.25]{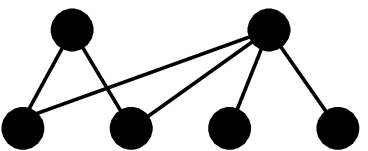} & 304 &\textbf{{\small 64461 }}& 21+18&  \includegraphics[scale=.25]{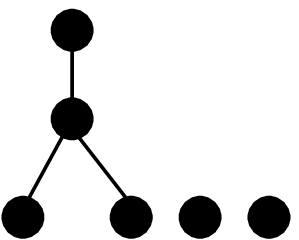} & 490 &\textbf{{\small 213428}}& 23+22\\
\hline
 \includegraphics[scale=.25]{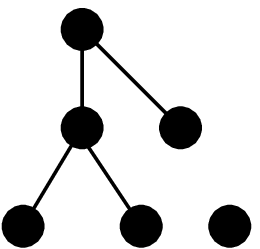} & 325 &\textbf{{\small 64004 }}& 20+18&  \includegraphics[scale=.25]{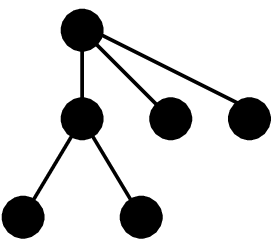} & 298 &\textbf{{\small 63640}}& 19+18\\
\hline
 \includegraphics[scale=.25]{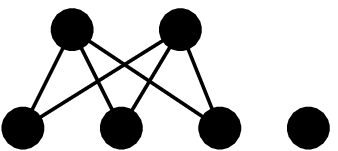} & 255 &\textbf{{\small 20984 }}& 20+15&  \includegraphics[scale=.25]{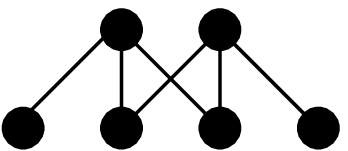} & 218 &\textbf{{\small 20392}}& 19+14\\
\hline
 \includegraphics[scale=.25]{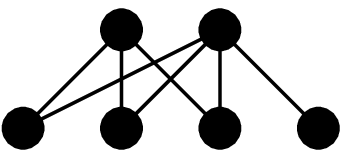} & 191 &\textbf{{\small 20184 }}& 18+14&  \includegraphics[scale=.25]{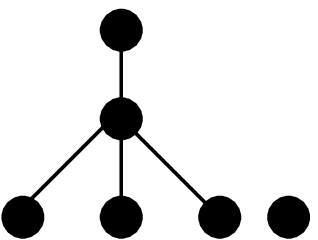} & 209 &\textbf{{\small 20379}}& 18+14\\
\hline
 \includegraphics[scale=.25]{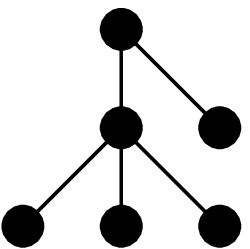} & 188 &\textbf{{\small 20181 }}& 17+14&  \includegraphics[scale=.25]{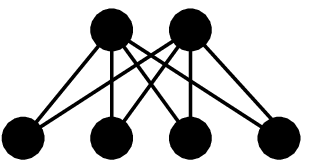} & 170 &\textbf{{\small 19986}}& 17+14\\
\hline
 \includegraphics[scale=.25]{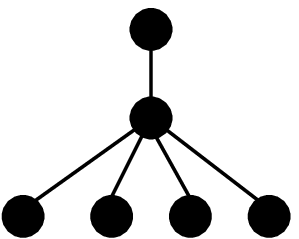} & 168 &\textbf{{\small 19984 }}& 16+14& \includegraphics[scale=.25]{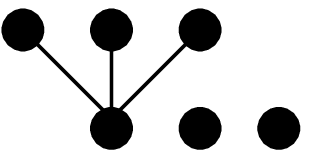} & 9944 & & 34+133\\
\hline
\includegraphics[scale=.25]{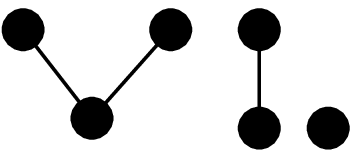} & 2024 &\textbf{{\small2610806855}}& 28+51& \includegraphics[scale=.25]{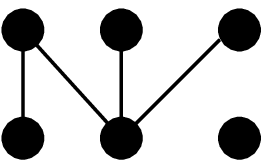} & 1195 &\textbf{{\small 179700889 }} & 26+43\\
\hline
 \includegraphics[scale=.25]{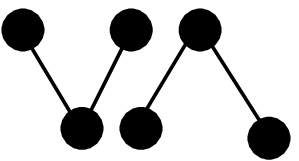} & 596 &\textbf{{\small 153926}} & 23+22&  \includegraphics[scale=.25]{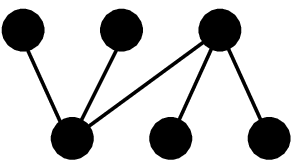} & 428 &\textbf{{\small 121130}} & 22+22\\
\hline
\includegraphics[scale=.25]{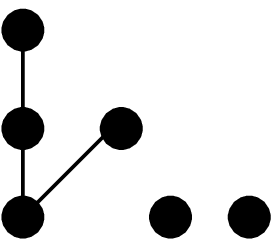} & 1326 &\textbf{{\small 296198143}} & 26+45&  \includegraphics[scale=.25]{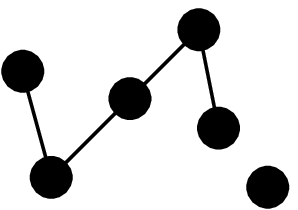} & 472 &\textbf{{\small 138454}} & 22+22\\
\hline
 \includegraphics[scale=.25]{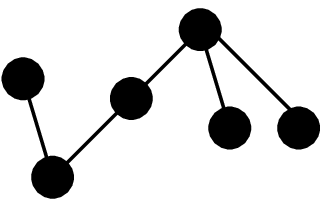} & 318 & \textbf{{\small 63872}} & 20+18& \includegraphics[scale=.25]{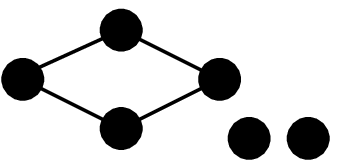} & 492 &\textbf{{\small 210044}} & 22+23\\
\hline
 \includegraphics[scale=.25]{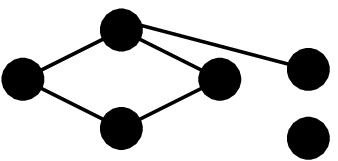} & 325 & \textbf{{\small 63943}} & 20+18&  \includegraphics[scale=.25]{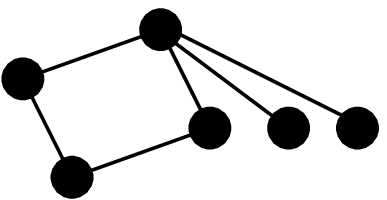} & 298 &\textbf{{\small 63640}} & 19+18\\
\hline
\includegraphics[scale=.25]{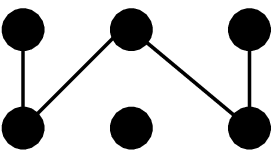} & 670 & \textbf{{\small 434366}} & 24+26& \includegraphics[scale=.25]{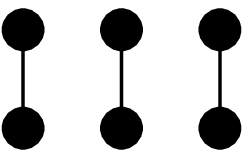} & 987 &\textbf{{\small 1007808}} & 25+27\\
\hline
 \includegraphics[scale=.25]{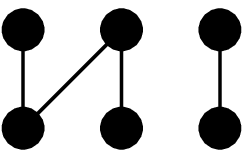} & 434 & \textbf{{\small 14616}} & 22+15& \includegraphics[scale=.25]{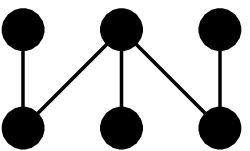} & 243 &\textbf{{\small 3311}} & 20+11\\
\hline
 \includegraphics[scale=.25]{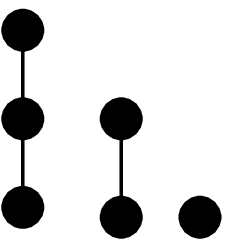} & 488 & \textbf{{\small 60962}} & 22+18& \includegraphics[scale=.25]{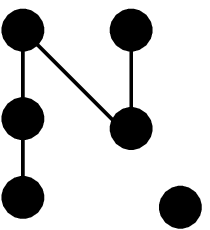} & 273 &\textbf{{\small 32449}} & 20+12\\
\hline
 \includegraphics[scale=.25]{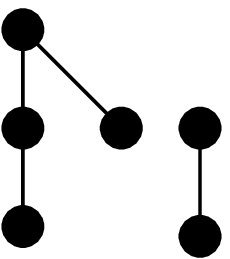} & 234 & \textbf{{\small 2895}} & 19+9&\includegraphics[scale=.25]{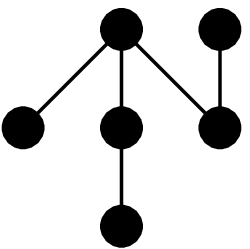} & 184 &\textbf{{\small 2626}} & 18+9\\
\hline
 \includegraphics[scale=.25]{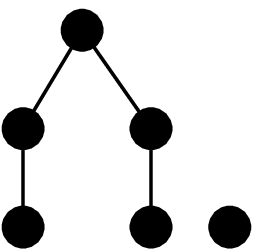} & 194 & \textbf{{\small 2665}} & 18+9& \includegraphics[scale=.25]{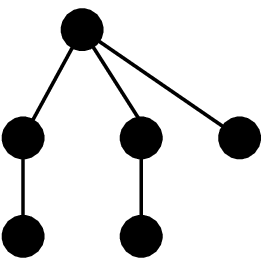} & 174 &\textbf{{\small 2604}} & 17+9\\
\hline
\includegraphics[scale=.25]{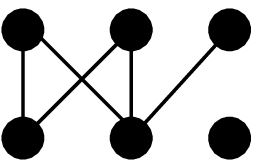} & 194 & \textbf{{\small 2665}} & 18+9& \includegraphics[scale=.25]{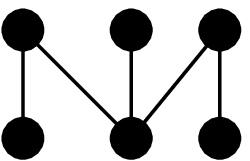} & 243 &\textbf{{\small 3311}} & 20+11\\
\hline
 \includegraphics[scale=.25]{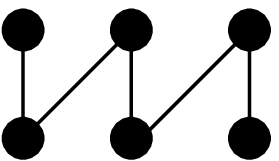} & 188 &  756 & 18+7& \includegraphics[scale=.25]{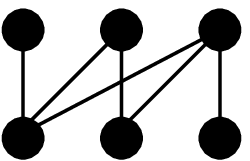} & 138 &584 & 18+7\\
\hline
 \includegraphics[scale=.25]{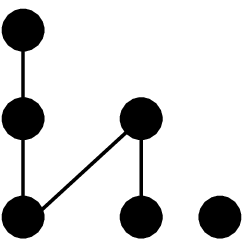} & 273 & \textbf{{\small 4936}} & 20+12& \includegraphics[scale=.25]{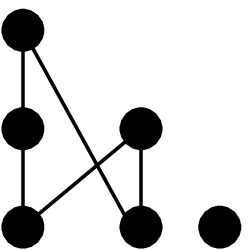} & 154 & 649 & 18+7\\
\hline
\includegraphics[scale=.25]{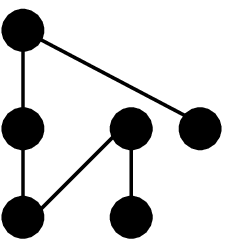} & 127 &  415 & 17+5& \includegraphics[scale=.25]{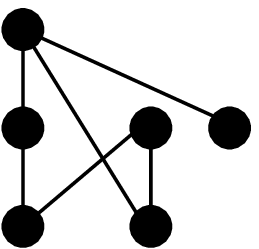} & 100 &361 & 16+5\\
\hline
\includegraphics[scale=.25]{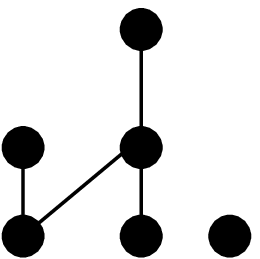} & 167 &  1060 & 18+8& \includegraphics[scale=.25]{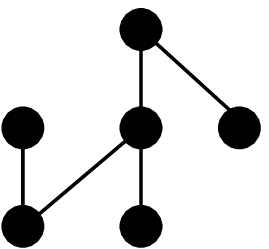} & 104 & 369 & 16+5\\
\hline
 \includegraphics[scale=.25]{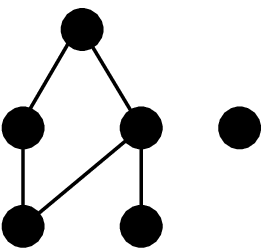} & 108 &  377 & 16+5& \includegraphics[scale=.25]{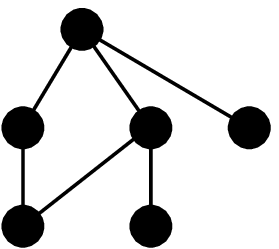} & 94 & 353 & 15+5\\
\hline
 \includegraphics[scale=.25]{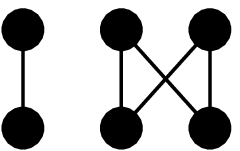} & 198 & 622 & 19+6& \includegraphics[scale=.25]{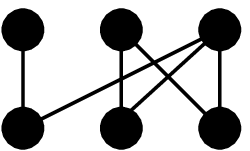} & 108 & 243 & 17+4\\
\hline
 \includegraphics[scale=.25]{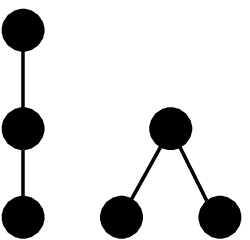} & 180 & 821 & 18+6& \includegraphics[scale=.25]{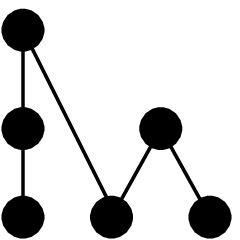} & 100 & 216 & 16+3\\
\hline
 \includegraphics[scale=.25]{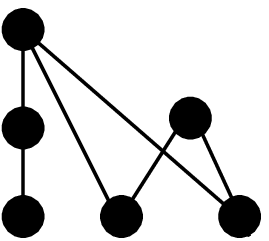} & 80 &  190 & 15+3& \includegraphics[scale=.25]{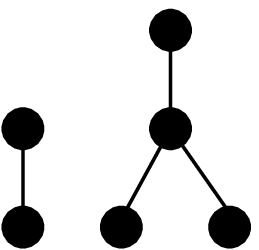} & 110 & 242 & 16+3\\
\hline
 \includegraphics[scale=.25]{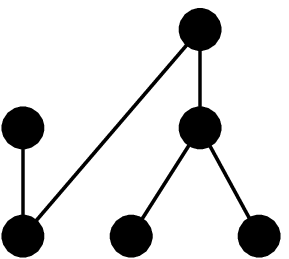} & 83 &  195 & 15+3& \includegraphics[scale=.25]{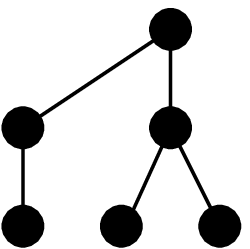} & 76 & 186 & 14+3\\
\hline
 \includegraphics[scale=.25]{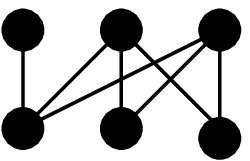} & 81 &  195 & 16+4& \includegraphics[scale=.25]{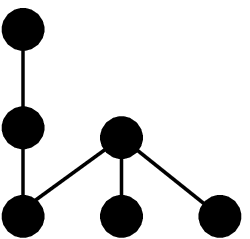} & 130& 686 & 17+6\\
\hline
\includegraphics[scale=.25]{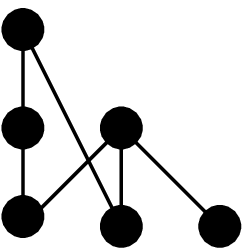} & 73 &  170 & 15+3& \includegraphics[scale=.25]{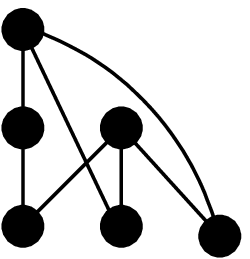} & 59& 151 & 14+3\\
\hline
 \includegraphics[scale=.25]{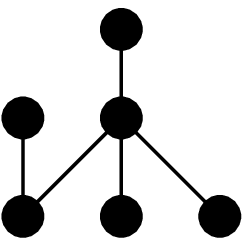} & 63 &  157 & 14+3& \includegraphics[scale=.25]{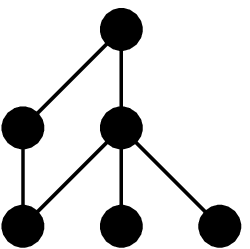} & 56& 148  & 13+3\\
\hline
 \includegraphics[scale=.25]{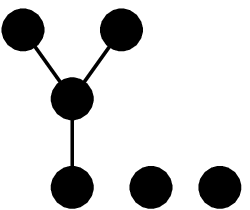} & 490 & \textbf{{\small 213428}} & 22+23& \includegraphics[scale=.25]{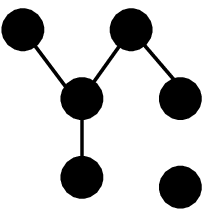} & 167& 1060  & 18+8\\
\hline
 \includegraphics[scale=.25]{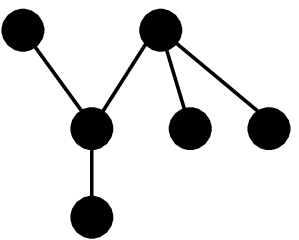} & 110 &  639 & 16+6& \includegraphics[scale=.25]{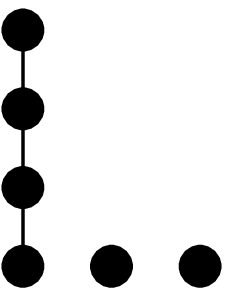} & 194& 2784 & 18+10\\
\hline
 \includegraphics[scale=.25]{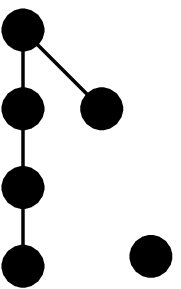} & 119 &  661 & 16+6&\includegraphics[scale=.25]{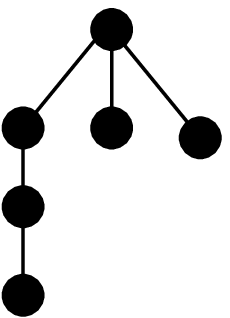} & 104& 630  & 15+6\\
\hline
 \includegraphics[scale=.25]{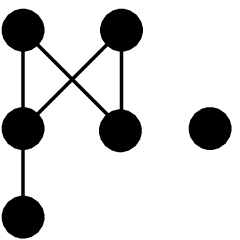} & 97 &  230 & 16+4& \includegraphics[scale=.25]{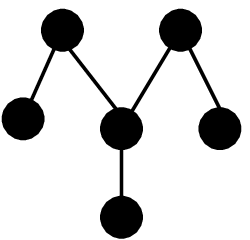} & 78& 178  & 15+3\\
\hline
\end{longtable}

Notice that the cardinality of $FM_3(P,\le)$ is printed boldface whenever $FM(P,\le)$ is infinite. As mentioned, otherwise the two cardinalities coincide. Thus e.g. $FM$\big( \includegraphics[scale=.25]{q57.eps}\big) $\cong$
$FM_3$\big( \includegraphics[scale=.25]{q57.eps}\big) has $756$ elements, $18$ factors $D_2$, $7$  factors $M_3$ and length $32$.

Call a poset $P$ {\it good} if  it does neither contain {\bf 1+1+1+1} nor {\bf 1+2+2} as subposet, and whence induces a finite lattice $FM(P,\le)$. All $1101$ good $7$-elements posets $P$ and their cardinalities  $|FM(P,\le)|$ are listed in [6, 8.2] as well.


\begin{thebibliography}{10}

\bibitem{berman}
J.~Berman and B.~Wolk,
\emph{Free lattices generated in some small varieties}, Algebra Universalis 10 (1980) 269-289.

\bibitem{Wi}
R.~Wille, \emph{On lattices freely generated by finite partially ordered sets},  Coll. Math. Soc. Janos Bolyai 1975, p.581-593.

\bibitem{Wi2}
R.~Wille, \emph{Subdirekte Produkte vollst\"{a}ndiger Verb\"{a}nde}, Journal fur die reine und angewandte Mathematik 283 (1976) 53-70.

\bibitem{Wi3}
R.~Wille, \emph{\"{U}ber modulare {V}erb\"{a}nde die von einer endlichen halgeordneten {M}enge frei erzeugt werden}, Math. Z. 131 (1973) 241-249.


\bibitem{W1}
M.~Wild, \emph{Computing various types of lattices freely generated by posets}, Note di matematica e fisica 10  (1999) 99-128.

\bibitem{yves}
J.Y.~Semegni, \emph{On the computation of freely generated modular lattices},  
Ph.D. thesis, University of Stellenbosch, 2008.

\bibitem{hw}
C.~Herrmann and M.~Wild \emph{Acyclic modular lattices and their representations}, Journal of Algebra, vol. 136 (1991) 365-396.

\bibitem{peter}
P.~ Luksch and M.~ Petkov\~scaronek, \emph{An explicit formula for $FM(1+1+n)$}, Order, 6(4) (1990) 319-324.

\bibitem{duquenne} 
V.~ Duquenne, \emph{Contextual implications between attributes and some representation properties for finite lattices}, Beitr\"{a}ge zur Begriffsanalyse, Wissenschaftsverlag Mannheim, Germany, 1987.

\bibitem{lux}
K.~ Lux and J.~ Muller and M.~ Ringe, \emph{Peakward condensation and submodule lattices: An application of the Meat-axe}, J. Symbolic Comput. 17 (1994) 529-544.

\bibitem{bartenschlager}
G.~Bartenschl\"{a}ger,
\emph{Free bounded distributive lattices generated by finite ordered sets}, PhD Tech. Hochschule Darmstadt 1994.

\bibitem{mw1}
M.~Wild,
\emph{Variations and applications of the implication $n$-algorithm}, in preparation.

\bibitem{mw2}
M.~Wild,
\emph{Optimal implicational bases for finite modular lattices}, Quaestiones Mathematicae 23 (2000) 153-161.

\bibitem{mw3}
M.~Wild,
\emph{Computing all maximum cardinality anticliques}, in preparation.
\end{thebibliography}
\end{document}